\documentclass[11pt]{article}
\usepackage{amssymb,latexsym,amsmath,amsthm,graphicx,subfigure}
\usepackage[percent]{overpic}

\makeatletter 
\@addtoreset{figure}{section}
\def\thefigure{\thesection.\@arabic\c@figure}
\def\fps@figure{h, t}
\@addtoreset{table}{bsection}
\def\thetable{\thesection.\@arabic\c@table}
\def\fps@table{h, t}
\@addtoreset{equation}{section}

\newif\ifamsfonts\amsfontstrue \ifamsfonts
\font\twlbbb=msbm10 scaled\magstep1 \font\egtbbb=msbm8
\font\sixbbb=msbm6
\newfam\mathbbfam
\textfont\mathbbfam=\twlbbb \scriptfont\mathbbfam=\egtbbb
\scriptscriptfont\mathbbfam=\sixbbb
\makeatother

\newtheorem{theorem}{Theorem}[section]
\newtheorem{proposition}[theorem]{Proposition}
\newtheorem{lemma}[theorem]{Lemma}
\newfont{\tenbi}{cmbxti10}

\newcommand{\todo}[1]{\vspace{5 mm}\par \noindent
\marginpar{\textsc{ToDo}}
\framebox{\begin{minipage}[c]{\textwidth}
\tt #1 \end{minipage}}\vspace{5 mm}\par}

\renewcommand{\epsilon}{\varepsilon}

\def\fder#1#2{{\frac {d #1}  {d #2}}}
\def\pder#1#2{{\frac {\partial #1} {\partial #2}}}

\def\coad#1#2{{\rm ad}^*_{#1}\, {#2}}
\def\Ad#1#2{{\rm Ad}_{#1}\, {#2}}
\def\Coad#1#2{{\rm Ad}^*_{#1}\, {#2}}

\begin{document}
\title{Discrete Nonholonomic LL Systems on Lie Groups
\footnote{AMS Subject Classification 37J60, 37J35, 70H45}} 
\author{Yuri N. Fedorov \\
 Department of Mathematics and Mechanics 
 \\ Moscow Lomonosov University, Moscow, 119 899, Russia \\ 
e-mail: fedorov@mech.math.msu.su \\ 
and \\ 
 Department de Matem\`atica I, \\ 
Universitat Politecnica de Catalunya, \\
Barcelona, E-08028 Spain \\
e-mail: Yuri.Fedorov@upc.es \\ 
and \\
Dmitry V. Zenkov
\\ 
Department of Mathematics \\
North Carolina State university \\ 
Raleigh, NC 27695 \\ 
e-mail: dvzenkov@unity.ncsu.edu} 

\date{\small September 21, 2004}
\maketitle 

\vspace{-2em}

\abstract{
\noindent
This paper applies the recently developed theory of discrete
nonholonomic mechanics to the study of discrete nonholonomic
left-invariant dynamics on Lie groups. The theory is illustrated with
the discrete versions of two classical nonholonomic systems, the
Suslov top and the Chaplygin sleigh. The preservation of the reduced
energy by the discrete flow is observed and the discrete
momentum conservation is discussed. 


}

\tableofcontents

\section{Introduction}
The theory of variational integrators for Lagrangian and Hamiltonian
systems originated in \cite{Ve1988}, \cite{Ve1991}, and \cite{MosVes}. It
was further developed by a number of authors (see e.g. \cite{BS1}, 
\cite{MCL1993} \cite{MarPekShk}, \cite{WeMa1997}, and \cite{MaWe2001}
for a more complete list of references and history). A very important
feature of variational integrators is the {\em discrete momentum
  preservation}: if the original continuous-time system has a symmetry
and conserves the momentum map, so does the associated discrete-time
mechanical system.

In \cite{CM1}, \cite{deLeon1} 
 the theory was extended to the  Lagrangian systems with nonholonomic
 constraints. 
In particular, it was shown in \cite{CM1} that the discrete-time
nonholonomic system conserves the {\em spatial} momentum in the case
of {\em horizontal symmetry} (see \cite{BKMM} for the definition of
the horizontal symmetry).  However, the case of horizontal symmetry is
not typical in nonholonomic mechanics. Apparently, Chaplygin
\cite{Ch1911} was the first to observe the link between symmetry and
conservation of the components of momentum relative to the {\em
  moving} frame (see also \cite{Z2003} and references
therein). Therefore, it is natural to ask whether the discrete
momentum is preserved by the discrete-time nonholonomic system
associated with a momentum-preserving continuous-time system. A
closely related property is the existence of an invariant measure. The
continuous-time nonholonomic systems generically are not
measure-preserving (see \cite{Koz1} and \cite{ZBl2002} for
details). The next version of this paper will address the
measure-preservation property for the discrete-time nonholonomic
systems.

The goal of this paper is to study the properties of the numerical
variational integrators for a nonholonomic mechanical system whose
configuration space is a Lie group $G$. Here we consider LL systems,
that is, we assume that both the Lagrangian 
and the constraint distribution 
are invariant with respect to the induced left action of $G$ on $TG$.

The paper is organized as follows: Section 2 gives a brief overview of
both continuous and discrete-time nonholonomic systems. 

In Section 3 we develop the theory of discrete left-invariant
nonholonomic systems on Lie groups $G$.  The fact that the constraints
on $G\times G$ are left-invariant enables us to reduce the dynamics on
a smooth {\it admissible displacement subvariety} ${\cal S}\subset G$,
which is chosen to be the exponent of a linear subspace ${\mathfrak
  d}$ of the Lie algebra $g$ of $G$.
Under the discrete Legendre transformation, ${\cal S}$ gives rise to a
discrete momentum locus ${\cal U}$ in the coalgebra $g^*$. In contrast
to continuous nonholonomic systems, the locus is not a linear subspace
in $g^*$, but rather a nonlinear subvariety.  The dynamics is then
described by the discrete Euler--Poincar\'e--Suslov equations
that generate a (generally multivalued) map from ${\cal U}$ onto
itself.

In Sections 4 and 5 we review the dynamics of the two classical
nonholonomic LL systems on the Lie groups $SO(3)$ and $SE(2)$, the
Suslov problem and the Chaplygin sleigh respectively, as well as their
multidimensional generalizations.

In Sections 6 and 7 we construct the discretizations of the above
problems as multi-valued maps on certain two-dimensional
non-orientable subvarieties of $SO(3)$ and $SE(2)$.  It is shown that
the discrete model retains such a distinct feature of the
continuous-time dynamics as the existence of heteroclinic trajectories
that connect the two one-parameter families of relative equilibria of
the system.  If, for special values of parameters, the continuous-time
system is momentum/measure preserving, then so is its discrete analog.

Moreover, it appears that in both discretizations the corresponding
reduced constrained energy is preserved as well.  This conservation
law replaces the momentum conservation in the general case and seems
to be quite unexpected, since generically the discrete variational
integrators do not preserve the energy and this property does not
change in the nonholonomic case.


\section{Lagrangian Mechanics with Nonholonomic Constraints}
In this section we briefly discuss the main concepts of nonholonomic
dynamics. For a complete exposition see \cite{Bloch2003} and \cite{BKMM}. 

\paragraph{The Euler--Lagrange Equations for Nonholonomic Systems.}
A nonholonomic Lagrangian system is a triple $(Q,L,{\mathcal D} )$,
where $Q$ is a smooth $n$-dimensional manifold called the {\em
  configuration space}, $ L:TQ \to \mathbb R $ is a smooth function
called the {\em Lagrangian}, and ${\mathcal D} \subset TQ$ is a $k$-dimensional
{\em constraint distribution}. Let $ q=(q ^1 ,\dots, q ^n)$ be local
coordinates on $Q$. In the induced coordinates $ (q, \dot q)$ on the
tangent bundle $ TQ $ we write $ L(q ,\dot q) $. It is assumed that
the Lagrangian is {\em hyperregular}, {\em i.e.,} the map
\[
\frac{ \partial  L}{\partial \dot q } : TQ \to T^* Q
\]
is invertible (see \cite{MaRa1999}).

A curve $q(t) \in Q$ is said to {\em satisfy the
constraints} if $\dot{q}(t)\in{\mathcal D}_{q(t)}$ for all~$t$.
The equations of motion are
given by the following La\-gran\-ge--d'Al\-em\-bert
principle:
\emph{The {La\-gran\-ge--d'Al\-em\-bert equations of
motion} for the system are those determined by
\begin{equation}\label{LD.eqn}
\delta \int^b_a L(q ^i , \dot{q} ^i ) \, d t = 0,
\end{equation}
where we choose variations $\delta q(t)$ of the curve $q(t)$ that
satisfy $\delta q(a) = \delta q(b) = 0$ and $\delta q(t) \in {\mathcal
  D}_{q(t)}$ for each $t$ where $ a\leq t\leq b$.}  This principle is
supplemented by the condition that the curve itself satisfies the
constraints. Note that we take the variation {\em before} imposing the
constraints; that is, we do not impose the constraints on the family
of curves defining the variation. This is well known to be important
to obtain the correct mechanical equations (see Bloch, Krishnaprasad,
Marsden, and Murray~\cite{BKMM} for a discussion and references).

Assuming that the constraint distribution is specified by a set of $ n -
k $ differential forms $ A^j(q)$, $ j = 1, \dots, n - k $,
\begin{equation} \label{00}
{\mathcal D} = \{\dot q \in TQ \mid \langle A ^j(q), \dot q \rangle =
0,\ j=1,\dots,s=n-k\},
\end{equation}
equation \eqref{LD.eqn} implies
\begin{equation}
\frac{d}{dt}\frac{\partial
   L}{\partial \dot q} -
 \frac{\partial L}{\partial q} 
= \sum_{j=1}^s \lambda_j A ^j(q).
\label{1.1}
\end{equation}
Equations \eqref{1.1} are
 called the {\em Euler--Lagrange equations with multipliers.} Coupled
 with \eqref{00}, they give a complete description of the dynamics of
 the system.
\begin{lemma}
Equations \eqref{1.1} conserve the energy
\begin{equation} \label{energy.eqn}
E = \left\langle \pder{L}{\dot q} , \dot q \right\rangle - L.
\end{equation} 
\end{lemma}

\begin{proof}
Differentiating \eqref{energy.eqn} along the flow \eqref{1.1}, one obtains
\begin{align*} 
\dot E &= \left\langle \fder{}{t} \pder{L}{\dot q} , \dot q \right\rangle + 
\left\langle \pder{L}{\dot q} , \ddot q \right\rangle - 
\left\langle \pder{L}{q} , \dot q \right\rangle - \left\langle
\pder{L}{\dot q}, \ddot q \right\rangle 
\\
&= 
\left\langle \fder{}{t} \pder{L}{\dot q} - \pder{L}{q}, \dot
q\right\rangle
= \left\langle \sum_{j = 1}^s \lambda _j A ^j (q), \dot q
\right\rangle
= 0,
\end{align*} 
since $ \langle A ^j (q), \dot q \rangle = 0 $, $ j = 1, \dots, s $. 
\end{proof}

\paragraph{The Euler--Poincar\'e--Suslov Equations.}
Let the configuration space be an $n$-dimensional connected Lie group
$G$ with local coordinates $g$. Let $\mathfrak{g}$ be the Lie algebra
of $G$, that is, the tangent space $ T _e G $ at the
identity element $ e \in G $ supplied with an antisymmetric bracket
operation $ [ \cdot\, , \cdot ] : \mathfrak g \times \mathfrak g \to
\mathfrak g $.


Define an {\em LL system} on $G$ as a Lagrangian system
$(G,L,{\mathcal D} )$ with a left-invariant Lagrangian $ L : TG \to
\mathbb R$ and a left-invariant (generally nonintegrable) distribution
${\mathcal D}$ on the tangent bundle $TG$.

The Lagrangian $L:TG\to \mathbb R $  is left-invariant if and
only if $ L (g, \dot g) $ depends on $ (g, \dot g )$ through the
combination $ \omega = g ^{-1} \dot g $, {\em i.e.,} there exists a
function $ l : \mathfrak g \to \mathbb R $ called the {\em reduced
Lagrangian} such that $ L(g, \dot g) = l(\omega) $.

A distribution $ {\mathcal D} \subset TG $ is left-invariant if and
only if there is a subspace $ \mathfrak d \subset \mathfrak g $ such
that ${\mathcal D} _g= TL _g \,\mathfrak{d} \subset T_g \,G$ for any $
g \in G $.  Let $\mathfrak g^*$ be the dual of the Lie algebra 
and $ a^j $, $ j = 1, \dots, s $, be
independent elements of $ \mathfrak g ^{*} $ associated with the
subspace $\mathfrak d $, {\em i.e.,}
\[
\mathfrak d = \{ \xi \in \mathfrak g \mid \langle a ^j , \xi \rangle =
0,\ j = 1, \dots, s \}.
\]
Then the left-invariant constraints can be written as
\begin{equation}\label{linv_const.eqn}
\langle a ^j, \omega \rangle = 0, \qquad j = 1, \dots, s,
\end{equation}
where $\omega = g ^{-1} \dot g = TL _{q ^{-1} } \dot g$ is the {\em
  body velocity operator}.

Define the {\em body momentum\/} $ p : \mathfrak g ^{*} \to \mathbb R $
by the formula $p = \partial {l}/\partial {\omega}$.
According to \cite{Koz2}, 
the {\em reduced dynamics} of an LL system $(G,L,{\mathcal D})$ 
is governed by the {\em Euler--Poincar\'e--Suslov equations}
\begin{equation} \label{EPS.eqn}
\dot p = \coad{\omega }{p} + \sum _{j=1} ^s \lambda ^j a _j 
\end{equation}
coupled with the constraint equations \eqref{linv_const.eqn}.
The dynamics of the group variables $g$ is obtained by solving 
 the {\em reconstruction equation}
\begin{equation}\label{rec.eqn}
\dot g = TL _g \,\omega.
\end{equation}

\begin{theorem}\label{reduced_energy.thm}
The Euler--Poincar\'e--Suslov equations conserve the {\bf \em reduced
  constrained energy}
\[
{\mathcal E}  = \left[ \langle p, \omega \rangle - l (\omega )\right]  _{\omega
  \in \mathfrak d}.
\]
\end{theorem}
\begin{proof}
To prove this statement, observe that the reduced energy, $ \langle p,
\omega \rangle - l (\omega ) $, equals the energy as the Lagrangian is
left-invariant. Since $ \omega \in \mathfrak d $ throughout the
motion, the reduced constrained energy equals the energy along the
trajectories of \eqref{EPS.eqn} and therefore is preserved.
\end{proof}

Let the reduced Lagrangian $l(\omega)$ be 
the quadratic form $l=\frac 12 \langle \omega, \mathbb I \omega
\rangle$, where $\mathbb I : {\mathfrak g} \to {\mathfrak
  g}^*$ is a symmetric
non-singular {\em inertia operator}. In this case $p=\mathbb I \omega$. 
Then the constraints
(\ref{linv_const.eqn}) imply that $p$ lies in the subspace
$$ \mathfrak d ^\perp =\{ \langle a ^j, \mathbb I ^{-1} p \rangle = 0, \; j = 1,
\dots, s\}\subset {\mathfrak g}^* .
$$ It is often convenient to choose a basis $ e _1, \dots, e _n $ in
the Lie algebra $\mathfrak g$ such that 
$ a ^j = e ^{n - j+1} $, $ j = 1, \dots, s $. In such a basis, the
reduced constrained energy becomes
\begin{equation}\label{reducedenergy.eqn}
\frac 12 \sum_{i, j = 1}^s \mathbb I _{ij} \omega ^i \omega ^j = 
\frac 12 \sum_{i, j =
  1}^s J ^{ij} p _i p _j.
\end{equation} 
Here and elsewhere, the quantities $ J ^{ij} $ represent the
components of the {\em inverse constrained inertia operator}~$ \mathbb I\, _{| \mathfrak d} $.

\paragraph{Remark.} In the absence of constraints, equations
\eqref{EPS.eqn} become the Euler--Poincar\'e equations, which conserve
the {\it spatial} momentum $J = \Coad {g}{p}$. In the presence of
nonholonomic constraints, neither the spatial, nor body momentum is
conserved generically. However, in some cases the body momentum is
preserved. 
(The conditions for the body momentum preservation can be seen in \cite{Z2003}).

\paragraph{Discrete Mechanical Systems with Nonholonomic Constraints.}
According to \cite{CM1}, a discrete nonholonomic mechanical system
on $Q$  is specified by
\begin{itemize}
\item[(i)]
a {\em discrete Lagrangian} $L_d : Q\times Q \to {\mathbb R}$;

\item[(ii)] an $(n-s)$-dimensional distribution $\cal D$ on $TQ$;

\item[(iii)] a discrete constraint manifold ${\mathcal D}_d\subset
  Q\times Q$, which has the same dimension as $\cal D$ and satisfies the
  condition $(q,q) \in \mathcal D_d$ for all $q \in Q$.
\end{itemize}
The dynamics is given by the following {\em discrete
  Lagrange--d'Alembert principle} (see \cite{CM1}), 
\[
\sum _{k = 0} ^{N-1}
\bigg(D _1 L _d ( q _k , q _{k+1}) + D _2 L _d (q _{k-1} , q _k)\bigg )\,
 \delta q_k = 0, \quad \delta q _k \in {\mathcal D} _{q _k }, \quad
(q_k , q_{k+1}) \in {\mathcal D}_d.
\]
Here $D_1 L _d $ and $ D_2 L _d $ denote the partial derivatives of
the discrete Lagrangian with respect to the first and the second
inputs, respectively.

The discrete constraint manifold is usually specified by the {\em discrete
constraint functions}
\begin{equation} \label{gen.constr}
{\cal F}_j(q_k, q_{k+1})=0, \qquad j=1,\dots, s,
\end{equation}
which impose the restriction $(q_k,q_{k +1}) \in {\mathcal D}_d$
on the solution sequence $\{(q_k, q _{k+1} )\}$.

\paragraph{Remark.}
If the discrete Lagrangian is obtained from a continuous one,
$L(q,\dot q)$, via a discretization mapping $\Psi : Q\times Q\to TQ$
defined in a neighborhood of the diagonal of $Q\times Q$, i.e.,
$L_d=L\circ \Psi$, then the variety ${\mathcal D}_d$ must be {\it
  consistent} with the continuous distribution ${\mathcal D}$:
${\mathcal D}_d$ is locally defined by the equations $\nu^j \circ \Psi
= 0$, $ j= 1, \dots, s$.
  We emphasize that the
discretization mapping is not unique and hence there are many ways to
define the discrete Lagrangian $ L _d $ and the discrete constraint
manifold ${\mathcal D} _d $ for a given nonholonomic system $(Q, L,
{\mathcal D})$.\footnote{An alternative approach to the discretization
  of nonholonomic systems based on a modification of canonical
  transformations was proposed in \cite{deLeon1}.}

The dynamics of a discrete nonholonomic system is  represented by
sequences $ \{(q _k, q _{k+1})\}$ that satisfy the {\em
  discrete Lagrange--d'Alembert equations with multipliers}
\begin{equation}\label{dis_LA.eqn}
D_1 L_d (q_k, q_{k+1})+ D_2 L_d (q_{k-1},q_k)
= \sum_{j=1}^s \lambda_{j} ^k A ^j (q_k), \quad
{\mathcal F}_j(q_k, q_{k+1}) =0.
\end{equation}
As in the continuous-time case, these equations are equivalent the
discrete Lagrange--d'Alembert principle.

\paragraph{Remark.}
According to \cite{CM1}, equations \eqref{dis_LA.eqn} introduce a 
well-defined mapping 
$(q_{k-1}, q_{k})\mapsto (q_{k}, q_{k+1})$, 
if the $(n+s)\times (n+s)$ matrix
$$ \begin{pmatrix}
D_1 D_2 L _d(q_k ,q_{k +1}) & A^1(q_k) & \cdots & A^s(q_k)  \\
D_2 {\cal F}_1 (q_k ,q_{k +1}) & 0  & \cdots & 0 \\
\vdots & \vdots &  & \vdots \\
D_2 {\cal F}_s (q_k ,q_{k +1}) & 0  & \cdots & 0
\end{pmatrix}
$$ is invertible for each $(q_k ,q_{k +1})$ in a neighborhood of the
diagonal of $Q\times Q$.

\section{Discrete Euler--Poincar\'e--Suslov Equations}

\paragraph{Continuous and Discrete Left-Invariant Lagrangians.}
Assume that the configuration space is a Lie group $G$ and denote the
local coordinates in $G$ by $g$.  Let the discrete Lagrangian $ L _d :
G \times G \to G $ be invariant with respect to the left diagonal
action of $G$ on $ G \times G $:
\[
L _d (g\, g_k, g\, g_{k+1})=  L _d (g_k, g_{k+1})
\]
for any $ g\in G$.

Define the {\em incremental displacement} by the formula $W_k=
g_k^{-1} g _{k+1}\in G$.  Since $ L_d $ is left-invariant, there
exists a function $ l _d : (G \times G)/G \cong G \to \mathbb R $
called the {\em reduced discrete Lagrangian} such that $L_d ( g_k , g
_{k+1}) =l_d (W_k) $.

According to \cite{MarPekShk}, for a given continuous left-invariant
Lagrangian $L(g, \dot g) = l(g ^{-1} \dot g)$ its discrete analog $l
_d$ can be chosen in form
\[
l _d = l ( (\log W _k ) / h),
\]
where $ \log : G \to \mathfrak g $ is the (local) inverse of the
exponential map $ \exp : \mathfrak g \to G $ and $h \in \mathbb R_{+}$ is the given time step.

For a matrix group $G$, one can approximate $ (\log W _k ) / h $ with \\
$W_k - I \equiv g _k ^{-1} ( g _{k+1} - g _k)$, so that
\begin{equation} \label{approx}
L_d (g_k,g_{k+1}) = l\left(g _k ^{-1} ( g _{k+1} - g _k)/h\right ). 
\end{equation}
This will be our default choice for the groups $SO(n)$ and $SE(n)$ considered
in the next sections.

Similarly to \cite{BS1, MarPekShk}, we define the {\em discrete body momentum}
$ P_k : G \times G \to \mathfrak g ^{*} $ by the formula
$$
P_k = L_{g_k}^* D_2 L_d (g_{k-1}, g_k) \equiv L_{W_k}^* D\, l_d (W_k), 
$$ where $L_{g_k}^*$ is the induced left action $L^* _{g _k } : T ^* G
\to \mathfrak g^*$.  Equivalently, $P_k$ is defined by any of the
conditions: for any $\omega\in{\mathfrak g}$,
\footnote{The definition of the
  discrete momentum \eqref{disc.mom.dfn} accepted in this paper computes
  $ p _k $ as a function of $ W _k $ whereas the standard definition used in many publications,
$$
\langle \omega, P_k\rangle = \frac d{d\varepsilon} L_d 
\left(g_{k-1}, g_{k} e^{\omega\varepsilon}\right )
\bigg |_{\varepsilon=0} .
$$
makes $ p _k $ a function of $W_{k-1}$.}
\begin{equation}
\label{disc.mom.dfn} 
\langle \omega, P_k\rangle 
= - \frac d{d\varepsilon} L_d \left(g_k e^{\omega\varepsilon}, g_{k+1} \right )
\bigg |_{\varepsilon=0} \quad
\mbox{or} \quad \langle \omega, P_k\rangle 
= -\frac d{d\varepsilon} l_d \left( e^{\omega\varepsilon} W_{k}\right )\bigg |_{\varepsilon=0}.  
\end{equation}

In the unconstrained case, this defines the {\em discrete Legendre transformation} 
${\mathcal L} : (g_k, W_k)\in G\times G \to (g_k,P_k)\in G \times {\mathfrak g}^*$. 
The mapping $ {\mathcal L} $ is uniquely
invertible in a neighborhood of the set $\{(g, P) \in G \times
\mathfrak g ^{*} \mid P=0 \}$, but it may fail to be globally invertible.

In the presence of generic discrete constraints (\ref{gen.constr}), the displacement
$W_k$ is restricted to the {\em admissible displacement subvariety} 
\[
 {\mathcal V}_k = \{W_k \in G \mid {\mathcal F}_j (g_k, g_k W_k)=0 ,
 \ j=1,\dots, s\}
\]
As a result, the discrete momentum $P_k$ is restricted to an
$(n-s)$-dimensional subvariety ${\mathcal U}_k \subset {\mathfrak g}^*$,
the image of ${\mathcal L} (g_k, {\mathcal V}_k)$ in ${\mathfrak g}^*$.
In case of generic discrete constraints, for different $k$ the subvarieties 
${\mathcal U}_k$ are different.  

\paragraph{Discrete Left-Invariant Constraints.}
If the continuous constraint distribution ${\mathcal D}$ is
left-invariant, it is natural to require that the
discrete constraint manifold $ {\mathcal D} _d $ is invariant
with respect to the left diagonal action of $G$ on $ G \times G $, that is,
$$
{\mathcal F}_j(g\, g _k, g\, g_{k+1} )={\mathcal F}_j(q_k, q_{k+1})
\quad \text{for any} \quad g\in G, \quad j=1,\dots, s.
$$
This implies that there exist functions 
$ f _j : G \to\mathbb R $, $j = 1, \dots, s $, such that 
$$
{\mathcal F}_j(q_k, q_{k+1}) = f_j (W_k).
$$ 
Consequently, ${\cal D}_d \subset G \times G$ is completely defined
by the admissible displacement subvariety
$$
{\cal S}= \{ f_ 1 (W)=0,\; \dots, \; f_s (W)=0\}\subset G,
$$
namely ${\cal D}_d =\{g, g h \}$, $g\in G, h\in {\cal S}$. 

The submanifold
${\mathcal S}$ should pass through the identity element $I$ in $G$, and the
tangent space  $T \mathcal S _I $ at the identity should be ``horizontal'',
{\em i.e.,} it should coincide with the linear subspace ${\mathfrak d}\subset
{\mathfrak g}$ generating the left-invariant distribution on $TG$.

The second property suggests that ${\mathcal S} = \{ W \in G \mid \log W \in
\mathfrak d \} $.

Equivalently, $\cal S$ can be chosen a union of all
one-parameter subgroups $G_\omega$ generated by all admissible vectors
$\hat\omega \in {\mathfrak d}$. In other words, one can set ${\cal S}=
\exp {\mathfrak d}$.  However, in case of generic $\mathfrak d$, the set
$\exp {\mathfrak d}$ is not a subvariety of $G$.

In this paper we concentrate on the important case when $G$ contains a
subgroup $H$ generated by subalgebra ${\mathfrak h}\subset {\mathfrak g}$ 
such that there is a decomposition ${\mathfrak g}= {\mathfrak h}
\oplus {\mathfrak d}$ and $({\mathfrak g}, {\mathfrak h})$ forms a symmetric
pair, that is
\begin{equation} \label{cartan}
[{\mathfrak h}, {\mathfrak h} ]\subset {\mathfrak h}, \quad
[{\mathfrak d}, {\mathfrak d} ]\subset {\mathfrak h}, \quad 
[{\mathfrak h}, {\mathfrak d} ]\subset {\mathfrak d} . 
\end{equation} 

\begin{proposition} \label{S-Cartan}
\textup{(\cite{Helgason})} Under the condition (\ref{cartan})
the set ${\cal S}= \exp {\mathfrak d}$ is a smooth submanifold of $G$, 
which is either homeomorphic to the symmetric space $G/H$ or is a factor of $G/H$ 
by a finite group action. 
\end{proposition}

Under conditions  (\ref{cartan}) the set $\exp {\mathfrak d}$ is known as
the Cartan model of the symmetric space $G/H$. 

Notice that 
the tangent bundle $T\,S$ is not a subset of the left-invariant distribution $D\subset TG$,
since the latter is not integrable.

Under the Legendre transformation $\cal L$, the discrete momentum $P_k$ 
is confined to the subvariety 
\[
{\mathcal U} = \{ p \in \mathfrak g ^{*} \mid p = L _{W} ^{*} l _d '
(W), W \in {\mathcal S} \} \subset g^*,
\]
which now does not depend on $k$.
It appears that in the examples considered below 
the map ${\cal S} \mapsto {\cal U}$ is uniquely invertible almost everywhere on $\cal U$. 

\paragraph{Discrete Euler--Poincer\'e--Suslov Equations.}
Assume that the discrete Lagrangian $ L _d : G \times G \to \mathbb R$,
the discrete constraint distribution $ {\mathcal D} _d $, and the 
constraint distribution $ {\mathcal D} $ are left-invariant with
respect to the left action of $G$ on $ G \times G$ and $TG$,
respectively. 

Define the {\em action sum} and the {\em reduced action sum} by the
formulae 
\[
S _d = \sum _{k=0} ^{N-1} L_d (g _k, g _{k+1}) \quad \text{and} \quad 
s _d = \sum _{k=0} ^{N-1} l_d (W_k),
\]
respectively and rewrite the nonholonomic constraints (\ref{00}) as
a set of vanishing one-forms $A ^j (g) \dot g = 0$.

Following \cite{CM1}, consider variation of $S_d$ assuming that the
variations $ \delta g _k $ satisfy the conditions 
$ A ^j  ( g _k) \delta g _k = 0 $, $ j = 1, \dots, s$ and $\delta g _0 = \delta g _N = 0 $.

For the left-invariant constraints given by (\ref{linv_const.eqn})
the admissible discrete variations are those $ \delta g _k\in TG _{g _k }$ 
that satisfy the conditions
\begin{equation} \label{**}
\langle a ^j , g _k ^{-1} \delta g _k \rangle = 0, \qquad j = 1,
\dots, s, \qquad k = 1, \dots, N-1.
\end{equation}
The following theorem extends the result of \cite{BS1, MarPekShk} to the
nonholonomic setting.

\begin{theorem} \label{discrete_EPS}
Let $ L _d : G \times G \to \mathbb R  $ be a left-invariant
Lagrangian, $ l _d : G \to \mathbb R $ be the reduced Lagrangian, and
$ {\mathcal D} $ and $ {\mathcal D} _d $ be the compatible constraint
distributions on $ TQ $ and $ Q \times Q $, respectively. 
Then,  following statements are equivalent:
\begin{itemize}
\item[(i)]
The sequence $ \{ (g _k, g _{k+1} ) \} _{k = 0} ^{N-1} $ is a critical
point of the action sum $ S _d : G ^{N+1} \to \mathbb R $ for
arbitrary constrained variations.
\item[(ii)]
The sequence $ \{ (g _k, g _{k+1} ) \} _{k = 0} ^{N-1} $ satisfies the
discrete Euler--Lagrange equations with multipliers (\ref{dis_LA.eqn}) with
$q$ replaced by $g$, that is,
\begin{equation} \label{dis_LA_1.eqn}
D _1 L _d ( g _k , g _{k+1}) + D _2 L _d (g _{k-1} , g _k) = \sum_{j=1} ^s 
\lambda _k ^j A _j (g _k)
\end{equation} 
which are coupled with the discrete constraint equations
${\cal F}_j(g_k, g_{k+1})= 0$.

\item[(iii)] The sequence $ \{W_k \} _{k=0} ^{N-1} $ is a
  critical point of the reduced action sum $ s _d : G ^{N-1} \to
  \mathbb R $ with respect to variations $ \delta W_k $, induced
  by the constrained variations $ \delta g _k $, and given by
\begin{equation} \label{var_W}
\delta W_k = W_k \left[
  g _{k+1} ^{-1} \delta g _{k+1} - \Ad{ W_k^{-1}
  }{g _k ^{-1} \delta g _k} \right].
\end{equation}
\item[(iv)]
The sequence $ \{ W_k \} _{k=0} ^{N-1} $ satisfies the equations
\begin{equation}\label{discEP.eqn}
l _d ' (W_{k-1}) TL _{W_{k-1} } - l _d ' (W_k) 
TL _{W_k } \Ad {W_k^{-1}}{} = \sum _{j=1} ^s \lambda _k ^j a _j
\end{equation} 
coupled with the discrete constraint equations 
$$
f_j (W_k) = 0, \qquad j = 1,\dots, s, \quad k = 1, \dots, N-1.
$$
\end{itemize}
\end{theorem}

Proof of the theorem is given in the end of the section. 

We now rewrite \eqref{discEP.eqn} in the form of discrete momentum equations.
For any $ \omega \in \mathfrak g $,
\[
l _d ' (W _{k-1} ) TL _{W _{k-1} } \omega =
\langle l _d ' (W _{k-1} ), TL _{W _{k-1} } \omega  \rangle =
\langle 
 TL ^{*} _{W _{k-1} } l _d ' (W _{k-1} ), \omega \rangle
\]
and similarly
\[
l _d ' (W _k) TL _{W _k } \Ad {W _k^{-1}}{} \omega = 
\langle 
l _d ' (W _k) , TL _{W _k } \Ad {W _k}{} \omega \rangle = 
\langle \Coad {W _k^{-1}}{} TL ^{*} _{W _k } l _d ' (W _k) ,\omega \rangle .
\]
Therefore, in view of the definition of the discrete momentum (\ref{disc.mom.dfn}), 
\eqref{discEP.eqn} becomes {\bf \em discrete Euler--Poincar\'e--Suslov equations}
\begin{equation}
\label{disc.momentum.eqn}
P_{k+1} - \Coad {W_k} {P_{k}} = \sum _{j = 1} ^s \lambda_{k+1} ^j a_j,
\end{equation} 
where $W_k$ is restricted to ${\cal S}$ and $p_k \in {\cal U}\subset {\mathfrak g}^*$.

The above equations extend the discrete Euler--Poincar\'e equations
obtained in \cite{BS1, MarPekShk} to the case when the discrete
left-invariant constraints are present.  Thus, they represent a
discrete analog of (\ref{EPS.eqn}) and define a map ${\cal B}\, : \:
{\cal U} \mapsto {\cal U}: P_k \to P_{k+1}$, which is generally
multi-valued.  Given $P_k$, one evaluates $P_{k+1}$ by
\begin{description}
\item{1.} Finding $W_k$ by inverting the Legendre transformation;

\item {2.} Calculating $\hat P_k= \Coad {W_k} {P_k}$;

\item {3.} Choosing $P_{k+1}$ as one of the points of intersection of
  the $(n-s)$-dimensional subvariety $\cal U$ with the linear space
  span$(a_1,\dots,a_s)$ passing through $\hat P_k$.
\end{description}

If the map is multivalued, one needs to make a choice of a branch of
${\mathcal B}$. One natural way of doing this is to start from a value
of $P_{k}$ whose norm is small and to select $P_{k+1}$ of the smallest
norm.
\medskip

\noindent{\it Proof of Theorem \ref{discrete_EPS}}.
We first prove the equivalence of (i) and (ii) following \cite{CM1}. 
Recall that the variations $ \delta g _k $ vanish at $ k = 0 $ and $ k= N $.
Computing the first variation of the discrete action sum
$ S _d $, we obtain
\begin{align*}
\delta S _d & = \delta \sum _{k=0} ^{N-1} L _d (g _k, g _{k+1})
\\
&= \sum _{k=0} ^{N-1} D _1 L _d (g _k, g _{k+1}) \delta g _k 
+ \sum _{k=0} ^{N-1} D _2 L _d (g _k, g _{k+1}) \delta g _{k+1}
\\
&= \sum _{k=1} ^{N-1} D _1 L _d (g _k, g _{k+1}) \delta g _k 
+ \sum _{k=1} ^{N-1} D _2 L _d (g _{k-1} , g _k ) \delta g _k 
\\
&= \sum _{k=1} ^{N-1} \left(D _1 L _d (g _k, g _{k+1}) + D _2 L _d (g
_{k-1} , g _k ) \right) \delta g _k \, .
\end{align*} 
Here the variations $ \delta g_k $ are not independent and satisfy the
conditions $A _j ( g _k ) \delta g _k = 0$. 
Therefore, $ \delta S _d = 0 $ if and only if (\ref{dis_LA_1.eqn}) is fulfilled. 

Next, we prove that (i) is equivalent to (iii). Notice that $ L _d = l
_d \circ \pi $, where $\pi : G \times G \to (G \times G)/G \cong G $
is given by $ (g _k, g _{k+1}) \mapsto g _k ^{-1} g _{k+1}
$. Therefore
\[
\delta s _d  = \delta S _d .
\]
The variation $\delta W_k$ is computed to be 
\begin{align*}
\delta W_k &= \delta (g _k ^{-1} g _{k+1} ) =
g _k ^{-1} \delta g _{k+1} +  \delta g _k ^{-1} g _{k+1} 
\\
&= g _k ^{-1} \delta g _{k+1} - g _k ^{-1} \delta g _k g _k ^{-1} g _{k+1} \\
&= (g _k ^{-1}  g _{k+1}) (g _{k+1} ^{-1} \delta g _{k+1}) 
- (g _k ^{-1}  g _{k+1}) (g _k ^{-1}  g _{k+1}) ^{-1} (g _k ^{-1}
\delta g _k ) (g _k ^{-1} g _{k+1})
\\
&= g _k ^{-1}  g _{k+1} \left[
  g _{k+1} ^{-1} \delta g _{k+1} - \Ad{\left(g _k ^{-1} g _{k+1} \right) ^{-1}
  }{\left(g _k ^{-1} \delta g _k \right)} \right] ,
\end{align*}
which yields (\ref{var_W}).

To prove the equivalence of (iii) and (iv), we use (\ref{var_W}) to compute
\begin{align*}
\delta s _d &= \delta \sum _{k=0} ^{N-1} l _d (W _k) 
=  \sum _{k=0} ^{N-1} \delta l _d (W _k) 
= \sum _{k=0} ^{N-1} l _d '(W _k) \delta W _k
\\
&= \sum _{k=0} ^{N-1} l _d '(g _k ^{-1} g _{k+1})
g _k ^{-1}  g _{k+1} \left[
  g _{k+1} ^{-1} \delta g _{k+1} - \Ad{\left(g _k ^{-1} g _{k+1} \right) ^{-1}
  }{\left(g _k ^{-1} \delta g _k \right)} \right]
\\
&=
\sum _{k = 1} ^N l _d ' \left( g _{k-1} ^{-1}  g _k \right) 
 \left(g _{k-1} ^{-1}  g _k\right) \left( g _k  ^{-1} \delta g _k \right)
\\
& \qquad 
- \sum _{k=0} ^{N-1} l _d ' \left(g _k ^{-1} g _{k+1} \right)
\left( \left(g _k ^{-1} g _{k+1} \right)
  \Ad{\left(g _k ^{-1} g _{k+1} \right) ^{-1}}
  {\left(g _k ^{-1} \delta g _k \right)} \right)
\\
&= \sum _{k = 1} ^{N-1} 
\left[
  l _d ' \left( g _{k-1} ^{-1}  g _k \right) TL _{ g _{k-1}
    ^{-1} g _k } \vphantom{\Ad {\left( g _k ^{-1} g _{k+1} \right) ^{-1} } {}}
\right.
\\
&
\qquad 
\left.
- l _d ' \left(g _k ^{-1} g _{k+1} \right) 
TL _{g _k ^{-1} g _{k+1} }
\Ad {\left( g _k ^{-1} g _{k+1} \right) ^{-1} } {}
\right] 
\left( g _k ^{-1} \delta g _k \right).
\end{align*}
Since the variations $ \delta g _k $ satisfy the conditions (\ref{**}), 
$ \delta s _d = 0 $ if and only if item (iv) holds.



\section{The Suslov Problem and its Multidimensional Generalizations}
The most natural example of LL systems is the nonholonomic Suslov
problem, which describes the motion of a rigid body about a fixed
point under the action of the following nonholonomic constraint: the
projection of the angular velocity vector $ \overrightarrow
{\omega}\in {\mathbb R}^3$ to a certain {\it fixed in the body} unit
vector $\overrightarrow \gamma$ equals zero:
\begin{equation}
(\overrightarrow  \omega, \overrightarrow \gamma)=0. \label{ep3.24}
\end{equation}
The configuration space of the problem is the group $SO(3)$. Under the
identification of Lie algebras $({\mathbb R}^3,\times)$ and $(so(3),
[\cdot , \cdot])$, $\overrightarrow \omega$ and $\overrightarrow
\gamma$ correspond to elements of $so(3)$ and the coalgebra $so^*(3)$
respectively.

 Let $ {\mathbb I}\, : \,{\mathbb R}^3\mapsto {\mathbb R}^3$ be the
 inertia tensor of the body.  Then the Lagrangian equals $L=\frac 12
 (\overrightarrow \omega, {\mathbb I} \overrightarrow \omega)$ and the
 momentum $p$ is represented by the vector $ \overrightarrow
 M=(M_1,M_2, M_3)^T={\mathbb I}\overrightarrow \omega$.  The left
 action of the group $SO(3)$ on $T SO(3)$ leaves the kinetic energy of
 the body and the constraint (\ref{ep3.24}) invariant.

For the Suslov problem the Euler--Poincar\'e--Suslov equations
(\ref{EPS.eqn}) on $so(3)$ become
\begin{equation}
\frac {d}{dt}({\mathbb I}\overrightarrow \omega)
={\mathbb I}\overrightarrow \omega\times \overrightarrow \omega+\lambda\overrightarrow \gamma,  \label{3.25}
\end{equation}
where $\times $ denotes the vector product in ${\mathbb R}^3$ and $\lambda$ is
the Lagrange multiplier. 
Differentiating (\ref{ep3.24}), we find 
$$
\lambda =-( {\mathbb I}\overrightarrow \omega\times \overrightarrow \omega, 
{\mathbb I}^{-1} \gamma)/ (\overrightarrow \gamma,{\mathbb I}^{-1}\overrightarrow \gamma). 
$$
Therefore, (\ref{3.25}) can be represented as
$$
{d\over dt}({\mathbb I}\overrightarrow \omega) = \frac 1{(\gamma,{\mathbb I}^{-1}\overrightarrow \gamma)} 
{\mathbb I}^{-1}\overrightarrow \gamma  \times (({\mathbb I}\overrightarrow \omega\times \overrightarrow \omega)\times \overrightarrow \gamma) ,
$$
which, in view of (\ref{ep3.24}), is equivalent to
\begin{equation}
{d\over dt}({\mathbb I}\overrightarrow \omega)
=({\mathbb I}\overrightarrow \omega,\overrightarrow \gamma)\, \overrightarrow \omega\times {\mathbb I}^{-1}\overrightarrow \gamma.
\label{ep3.26}
\end{equation}
The Suslov system possesses the energy integral
\begin{equation}
( \overrightarrow \omega,{\mathbb I} \overrightarrow \omega) \equiv 
(\overrightarrow M,{\mathbb I}^{-1}\overrightarrow M )=
h ,\qquad h={\rm const} \label{ep3.27}
\end{equation}
and, as seen from (\ref{ep3.26}), it has a line of equilibria positions
$$
E=\{(\overrightarrow \omega,\overrightarrow \gamma)=0\}\cap \{({\mathbb I}\overrightarrow \omega,\overrightarrow \gamma)=0\}.  
$$
Note that in the principal basis, where 
${\mathbb I}=\mbox{diag }({\mathbb I}_1,{\mathbb I}_2,{\mathbb I}_3)$, the system
has the integral given by degenerate quadratic form 
\begin{equation} \label{cont_sus_int}
(\overrightarrow M, \hat{\mathbb I} \overrightarrow M), \quad
\hat{\mathbb I} = \begin{pmatrix}
{\mathbb I}_{2}\gamma _{3}^{2}+{\mathbb I}_{3}\gamma _{2}^{2} & 
-{\mathbb I}_{3}\gamma _{1}\gamma _{2} & 
-{\mathbb I}_{2}\gamma _{1}\gamma _{3} \\ 
-{\mathbb I}_{3}\gamma _{1}\gamma _{2} & 
{\mathbb I}_{1}\gamma _{3}^{2}+{\mathbb I}_{3}\gamma _{1}^{2} & 
-{\mathbb I}_{1}\gamma _{2}\gamma _{3} \\ 
-{\mathbb I}_{2}\gamma _{1}\gamma _{3} 
& -{\mathbb I}_{1}\gamma _{2}\gamma _{3} & {\mathbb I}_{1}\gamma
_{2}^{2}+{\mathbb I}_{2}\gamma_{1}^{2}
\end{pmatrix},
\end{equation}
which coincides with the restriction of (\ref{ep3.27}) onto the constraint plane 
$(\overrightarrow M, {\mathbb I}^{-1}\overrightarrow \gamma)=0$. 

In the basis where only one of the components of $\overrightarrow
\gamma$ is nonzero, say $\overrightarrow \gamma=(0,0,1)^T$, and the
inertia tensor is unbalanced, the integral (\ref{ep3.27}) can be
replaced by the reduced constrained energy integral
\begin{equation} \label{sus_int_1}
 {\mathbb I}_{22} M_1^2 -2 {\mathbb I}_{12} M_1 M_2 + {\mathbb I}_{11} M_2^2 .
\end{equation}

\begin{figure}[h,t]
\begin{center}
\includegraphics[width=.36\textwidth]{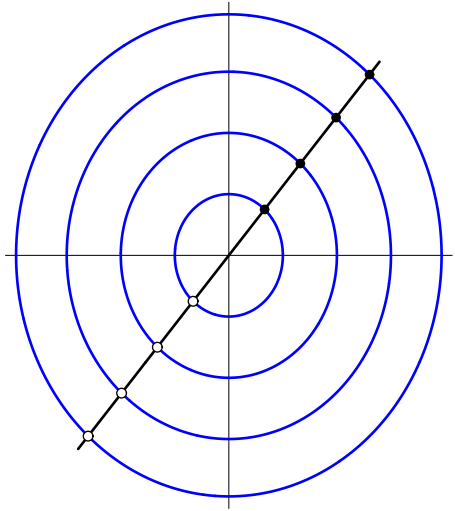}
\caption{\footnotesize The Momentum Dynamics for the Suslov Problem.}
\label{Suslov.fig}
\end{center}
\end{figure}

The dynamics of the two independent momentum components, $M_1$ and $M_2$, is illustrated in the Figure \ref{Suslov.fig}.
Because of the conservation law \eqref{sus_int_1}, 
the trajectories are the elliptic
arches that form the heteroclinic connections between the
asymptotically stable (filled dots) and unstable (empty dots)
equilibria.

As a result, the motion of the rigid body is
the asymptotic evolution from a permanent rotation about an axis fixed in the body frame 
to a permanent rotation about the same axis and with the
same angular velocity, but in the opposite direction. Note that {\it in space} 
the axes of the limit permanent rotations are different. 


The Suslov problem admits some natural multidimensional generalizations
studied in \cite{FeKo2,Jo2,Bl_Z}. 
The configuration space of an $n$-dimensional rigid body with a fixed point
 is the Lie group $SO(n)$. 
For a path $R(t)\in SO(n)$, the angular velocity of the body is
defined as the left-trivialization $\omega(t)=g^{-1}\cdot g(t) \in so(n)$.
 
The left-invariant metric on $SO(n)$ is given by non-degenerate inertia operator
${\mathbb I}\,:\, so(n)\to so(n)$. Then the Lagrangian of the free motion of
the body reads
\begin{equation} \label{cont.so-n.lagr}
L=\frac12\langle {\mathbb I}\omega,\omega\rangle,
\end{equation}
where now $\langle\cdot,\cdot\rangle$ denotes the Killing
metric on $so(n)$, $\langle X,Y\rangle=-\frac12\mbox{tr }(XY)$, $X,Y\in so(n)$.
For a ``physical'' rigid body, $\mathcal I\omega$ has the form
$J\omega+\omega J$, where $J$ is a symmetric $n\times n$ matrix called {\it mass
tensor} (see \cite{FeKo2}). 

Let $e_1,\dots,e_n$ be
the  orthogonal frame of unit vectors fixed in the body. 
What form may have a multi-dimensional analog of the condition (\ref{ep3.24})?
To answer this question, note that, instead of rotations {\it about an axis}
in the classical mechanics, in the $n$-dimensional case we have
infinitesimal rotations in the two-dimensional planes spanned by the
basis vectors $e_{i},e_{j}$,  $i,j=1,\ldots,n.$ 

Suppose, without loss of generality, that $\overrightarrow \gamma=(0,0,1)$ in (\ref{ep3.24}).  
Then this condition can be redefined as follows: only infinitesimal rotations 
in planes $(e_{1}, e_{3})$ and $(e_{2}, e_{3})$ are allowed. Hence,
it is natural to define the $n$-dimensional analog of Suslov's condition
in the following way: only infinitesimal rotations in the planes
$(e_{1}, e_{n}), \dots, (e_{n-1}, e_{n})$ (i.e., in the planes containing 
the vector $e_{n})$ are allowed. Thus, in the above basis, the angular velocity matrix
in the body must have the form
\begin{equation}
\omega =\begin{pmatrix} 
0            & \ldots & 0 & \omega_{1n} \\
\vdots       &        &   &  \vdots \\
     0       &        &   &  \omega_{n-1,n} \\
-\omega_{1n} & \ldots & -\omega_{n-1,n} & 0 \end{pmatrix} \, .
\label{matrix_constraint}
\end{equation}

This implies the constraints
\begin{equation}
\langle \omega, e_i\wedge e_j\rangle \equiv (e_i, \omega e_j)=0, \quad 1\le i<j\le n-1.
\label{Suslov_constraints}
\end{equation}
 As a result, the multidimensional Suslov problem is described by the EPS equations
on the Lie algebra $so(n)$
\begin{equation}
\frac{d}{dt} \left({\mathbb I}\omega\right) =[\mathcal I\omega,\omega ]
+\sum_{1<p<q\le n-1} \lambda_{pq} \, e_p\wedge e_q ,
\label{Suslov_eq}
\end{equation}
where the multipliers $\lambda_{pq}$ can be found by differentiating the constraints  
(\ref{Suslov_constraints}). 

Integrability of the system (\ref{Suslov_eq}),
(\ref{Suslov_constraints}) was proved, and its geometric properties
were studied in \cite{FeKo2}, whereas the reconstructed motion on the
group $SO(n)$ was described in \cite{Bl_Z}.


\section{Chaplygin Sleigh} 
Another example of a mechanical system governed by the
Euler--Poincar\'e--Suslov equations is the so-called Chaplygin sleigh,
the system introduced and studied in 1911 by Chaplygin \cite{Ch1911} (the work had
been actually finished in 1906, see also \cite{NeFu}).

\begin{figure}[h,t]
\begin{center}
\begin{overpic}[width=.5\textwidth]{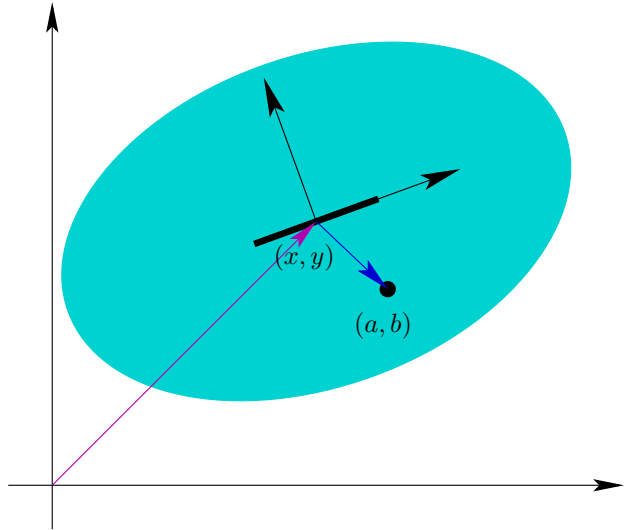}
\put(43,43){\small$(x,y)$}
\put(56,32){\small$(a,b)$}
\end{overpic}
\caption{\footnotesize The Chaplygin Sleigh}
\label{Chsleigh.fig}
\end{center}
\end{figure}

The sleigh is a rigid body moving on a horizontal plane supported at three
points, two of which slide freely without friction while the third is
a knife edge which allows no motion orthogonal to its direction, as shown in  
Figure \ref{Chsleigh.fig}. 

The configuration space of this dynamical system is the group of
Euclidean motions of the two-dimensional plane $ \mathbb R ^2 $,
$SE(2)$, which we parameterize with coordinates $ (\theta, x, y) $. As
the figure indicates, $\theta$ and $ (x,y) $ are the angular
orientation of the blade and position of the contact point of the
blade on the plane, respectively.

\paragraph{The Lagrangian and Constraint in the Body Frame.}
Introduce a coordinate system called the {\em body frame} by placing
the origin at the contact point and choosing the first coordinate axis
in the direction of the knife edge. Denote the angular velocity of the
body by $\omega =\dot\theta$, and the components of the linear velocity of the
contact point relative to the body frame by $v_1, v_2$. 
The set $ (\omega, v_1, v_2)$ is
regarded as an element of the Lie algebra $se(2)$. 

 The position of
the center of mass is specified by the coordinates $ (a, b) $ relative
to the body frame (we not assume here that the center of
mass lie along the blade direction as in some models). 
 We will see that $a$ is crucial to qualitative
behavior of the system while $b$ is irrelevant.

The Lagrangian equals the kinetic energy of the body, which is a sum
of the kinetic energy of the center of mass and the kinetic energy due
to the rotation of the body. Let $m$ and $J$ denote
the mass and moment of inertia of the sleigh relative to the contact
point. 
The position of the center of mass relative to the fixed (inertial) frame is
\[
(x + a \cos \theta - b \sin \theta, y + a \sin\theta + b \cos
\theta).
\]
Thus, the kinetic energy of the center of mass has the form 
\[
m \left[
\dot x ^2 + \dot y ^2 + (a ^2 + b ^2) \dot \theta ^2 + 2\dot \theta 
  \left(
    a (-\dot x \sin \theta + \dot y \cos \theta)
    - b (\dot x \cos \theta + \dot y \sin \theta)
\right) \right],
\]
or, using the body components of the angular and linear velocity,
\[
m \left[
(a ^2 + b ^2) \omega^2 + v_1^2 + v_2^2 - 2b
\omega v_1^2 + 2a \omega v_2 \right] .
\]
As a result,  the (reduced) Lagrangian is 
\begin{equation}
l = \frac12 \left[
 (J + m(a ^2 + b ^2)) \omega^2 + m (v_1^2 + v_2^2) 
- 2mb \omega v_1^2 + 2ma \omega v_2 \right].
\label{Caplag}
\end{equation}

Next, the constraint written relative to the body frame is $v_2= 0$.
Both the Lagrangian and constraint are invariant with respect to
the left action of $ SE(2)$ on $ TSE(2) $ as they depend on $ (g, \dot g)
$ through the combination 
\begin{equation}\label{se2_reconstruction.eqn} 
\Omega = g ^{-1} \dot g.
\end{equation} 

\paragraph{The Dynamics of Chaplygin Sleigh.}
In view of (\ref{Caplag}), the components of the body momentum are 
\begin{gather*}
p _\theta \equiv \pder{l}{\omega}=(J + m(a ^2 + b ^2)) \omega + 2m(a v_2-b v_1) , \\
p_1 \equiv \pder{l}{v_1} =m ( v_1-b \omega), \quad 
p_2 \equiv \pder{l}{v_2} =m ( v_2 +a \omega). 
\end{gather*}
The reduced dynamics of the Chaplygin sleigh is governed by the equations
\begin{equation} \label{EPS.se2.eqn}
\dot p _\theta = p _1 v_2 - p_2 v_1, \quad 
\dot p_1 = p_2 \omega, \quad 
\dot p _2 = - p_1 \omega + \lambda, 
\end{equation} 
which are the Euler--Poincar\'e--Suslov equations \eqref{EPS.eqn} on
the algebra $ se(2) $ coupled with the constraint $v_2= 0$.
Eliminating the variables $\Omega$ and the Lagrange multiplier
$\lambda$ from \eqref{EPS.se2.eqn}, one obtains the reduced dynamics
of the Chaplygin sleigh in the form of the momentum equation
\begin{equation}
\label{ch_sl_1.eqn}
\begin{aligned}   \dot p_\theta &= - \frac {a}{(J + ma^2)^2} 
  \left( p_\theta + b\, p_1 \right) 
  \left( m b \, p_\theta + (J + m(a^2 + b ^2)) \, p_1 \right), \\
  \dot p_1 &= \frac {ma}{(J + ma^2)^2}  \left( p_\theta + b \, p_1 \right) ^2 ,
\end{aligned}
\end{equation} 
which has the constrained energy integral
\begin{equation} \label{cont_sl_int}
 m p_\theta^2 +2 bm p_\theta p_1 + (J + m(a^2 + b ^2)) p_1^2 .
\end{equation}
In the case $b = 0 $ equations \eqref{ch_sl_1.eqn} become 
\begin{equation} \label{sleqns}
\dot p _\theta = - \frac{a \,p _\theta p _1 }{J + ma ^2 },
\qquad
\dot p _1 = \frac{ ma\, p _\theta ^2 }{(J + ma^2)^2}.
\end{equation}

We emphasize that the phase portrait of \eqref{ch_sl_1.eqn} 
is identical to that in the Suslov problem.
Indeed, if $a = 0$, the nonholonomic momentum $(p_\theta, p _1) $ is
conserved. Therefore, the body angular velocity $\omega$ and the
component of the body linear velocity along the blade $v_1$
are constants. The evolution of the configuration variables $ (\theta,
x, y) $ is determined from the {\em reconstruction equation}
\eqref{se2_reconstruction.eqn}, which reads
\begin{equation}\label{xyth.eqn}
\dot \theta =\omega, \qquad \dot x \cos \theta + \dot y \sin
\theta = v_1, \qquad - \dot x \sin \theta + \dot y \cos \theta =0. 
\end{equation} 
The solutions of \eqref{xyth.eqn} are 
\[
\theta = \theta_0 + \omega t, \quad x = x _0 
+ \frac{v_1}{\omega} \sin (\theta_0 + \omega t), \quad 
y = y_0 - \frac{v_1}{\omega} \cos (\theta _0 + \omega t) \quad \text{if} \quad
\omega \neq 0  
\]
and 
\[
\theta = \theta _0, \quad x = x _0 + v_1\cos\theta _0 \, t,
\quad y = y _0 + v_1 \sin \theta _0 \quad \text{if} \quad \omega = 0 .
\]
Therefore, the contact point of the blade and the plane
generically moves along a circle at a uniform rate.

\begin{figure}[h]
\begin{center}
\subfigure[$a \neq 0, b \neq 0$]{
\includegraphics[height=0.25\textwidth]{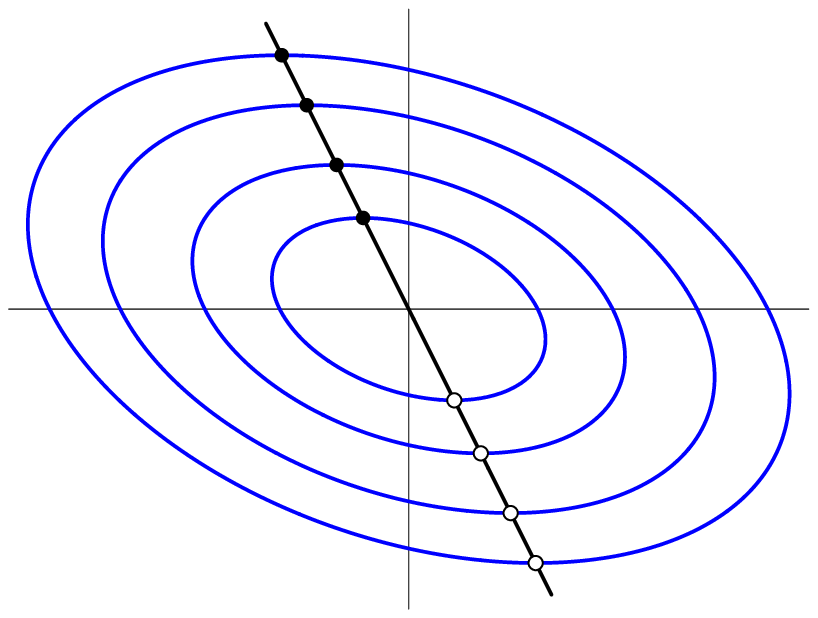}} \quad 
\subfigure[$a \neq 0, b = 0$]{
\includegraphics[height=0.25\textwidth]{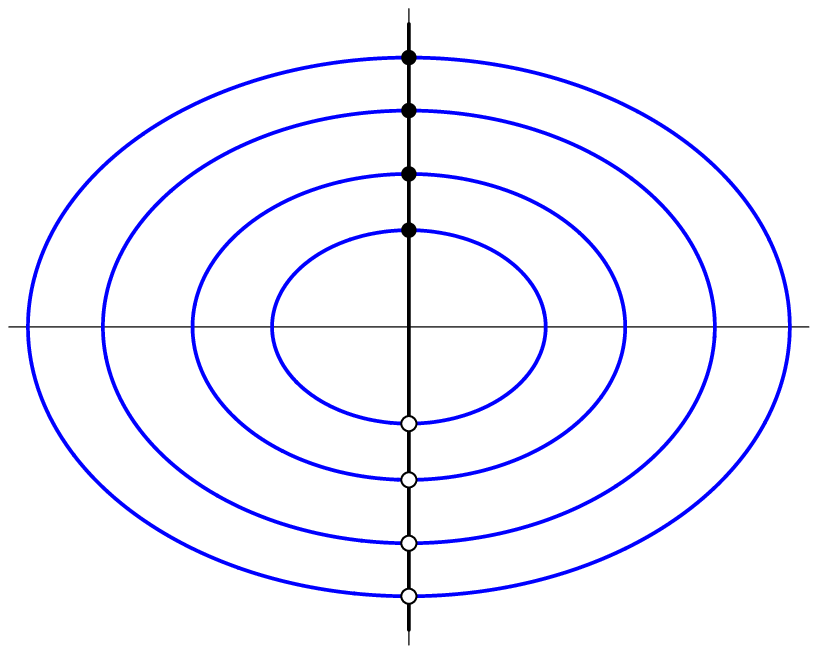}}
\end{center}
\caption{Momentum dynamics.}\label{mom_dynamics.fig}
\end{figure}

If $ a \neq 0 $, the dynamics \eqref{ch_sl_1.eqn} 
is integrable as the reduced energy is conserved.
The trajectories of \eqref{ch_sl_1.eqn}
 are either equilibria situated on the line 
$p_\theta+ b p_1=0$, or elliptic arches.\footnote{This follows from
  matching the trajectories and 
  the level curves of the reduced energy, which is a positive-definite
  quadratic form.}  
The equilibria located in the upper half plane are
asymptotically stable (filled dots in Figure \ref{mom_dynamics.fig})
whereas the equilibria in the lower half plane are unstable (empty
dots in Figure~\ref{mom_dynamics.fig}). The elliptic arches form
heteroclinic connections between the pairs of equilibria as shown in
Figure \ref{mom_dynamics.fig}.

\begin{figure}[h,t]
\begin{center}
\subfigure[$a \neq 0, b \neq 0$]{
\includegraphics[height=0.25\textwidth]{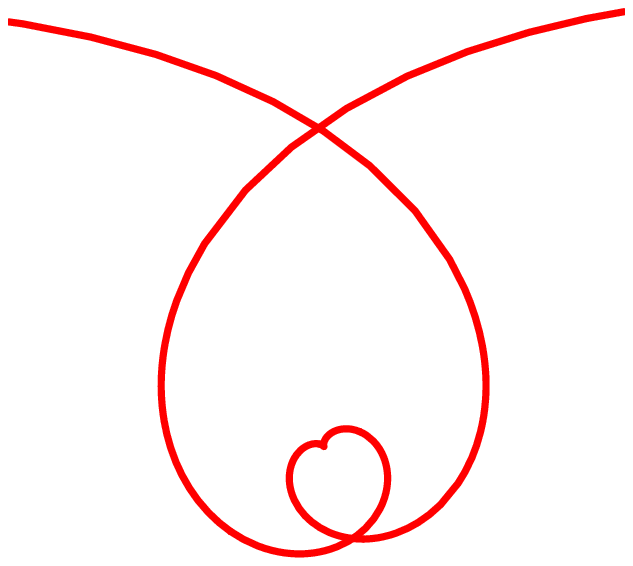}} \qquad 
\subfigure[$a \neq 0, b = 0$]{
\includegraphics[height=0.25\textwidth]{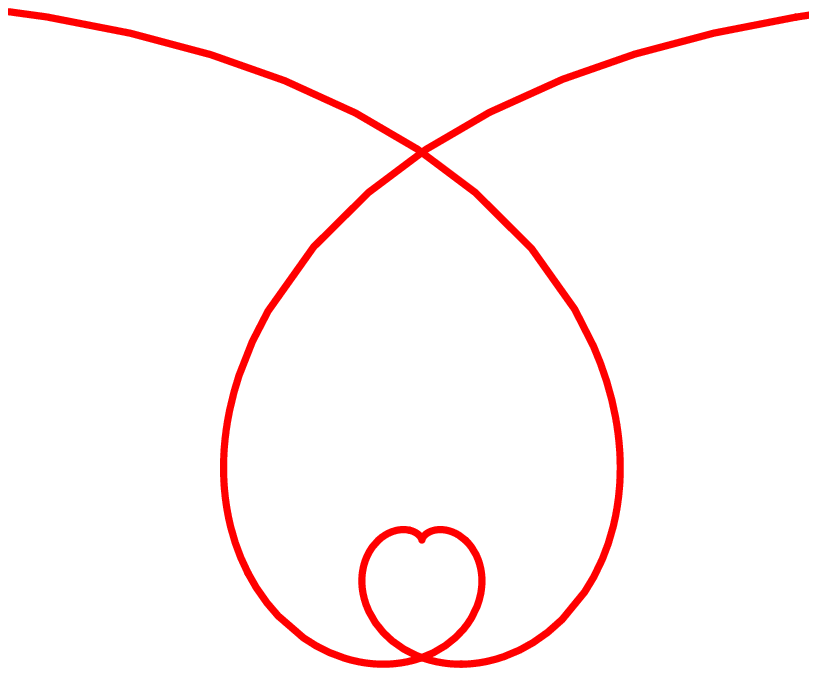}}
\end{center}
\caption{Generic trajectory of the blade.}
\label{generic_xy.fig}
\end{figure}

A generic trajectory of the contact point of the blade and the plane
has a cusp point (see Figure \ref{generic_xy.fig}). At the cusp, the speed of the
contact point, $ |v_1 | $, momentarily vanishes as the momentum
trajectory intersects the line $ mb p _1 + (J + m ( a ^2 + b ^2)) p _2= 0 $.

Since the group $SE(2)$ is 
a ``non-compact'' version of the group $SO(3)$, the dynamics of the
Chaplygin sleigh can be interpreted as a ``non-compact limit'' of the
dynamics of the Suslov problem.  Recall that the any non-equilibrium
trajectory of the Suslov top has a steady-state rotation as its
asymptotic dynamics. In a similar manner, a non-equilibrium state of
the Chaplygin sleigh asymptotically approaches a uniform straight-line
motions as $ t \to \pm \infty $. 

The shape of the generic trajectory of the contact point is
predetermined by the inertia of the body and the position of the
center of mass relative to the blade, and is independent of the
initial conditions. While the dynamics of the group variables $
(\theta, x, y) $ cannot be explicitly written, it is possible to
compute the angle between the asymptotic directions of the dynamics of
the contact point. See \cite{Ch1911} and \cite{NeFu} for details.

\paragraph{Multidimensional Chaplygin Sleigh.}
We now briefly discuss the generalized Chaplygin sleigh, which is an
$n$-dimensional rigid body moving in ${\mathbb R}^n$ in the presence
of certain nonholonomic constraints.

The configuration space of this dynamical system is the group $SE(n)$,
which has the structure of a semidirect product, $SE(n) = SO(n)
\,\circledS\, \mathbb R ^n$, so the group elements are written as
$(R,x)$, where $R \in SO(n)$ is the orthogonal rotation matrix of the
body and $x \in \mathbb R ^n$ is the position vector of its origin
$A$.  It is often convenient to represent the elements of $SE(n)$ by
means of $(n+1)\times (n+1)$ matrices of the form
$$
g(R, x) =
\begin{pmatrix}
R & x \\
0 & 1 \end{pmatrix},
$$ and the group operations for $SE(n)$ correspond to operations with
the matrices: the product of two such matrices corresponds to the
superposition of two Euclidean motions represented by these matrices
and the inverse matrix correspond to the inverse Euclidean motion.

The Lie algebra $se(n)$ of the group $SE(n)$ is the semidirect
product $so(n) \, \circledS \, \mathbb R ^n$ and it 
is isomorphic to the set of $(n+1)\times (n+1)$ matrices 
\[
\eta  = 
\begin{pmatrix}
\xi & v \\
0 & 0 
\end{pmatrix}, \quad \xi \in so(n), \quad v\in {\mathbb R}^n .
\]
The elements of $ se (n)$ are written as $ (\xi, v)$. 
The Lie bracket $[ \eta _1, \eta _2]$ in $se(n)$ is  $ \eta _1 \eta _2
- \eta _2 \eta _1 $, which yields 
\[
[(\xi _1, v _1), (\xi _2, v _2)] = ( [\xi _1, \xi _2], \xi _1 v _2 -
\xi _2 v _1 ).
\]

For a trajectory $g(t)\subset SE(n)$, the body velocity operator is
defined as the left-trivialization $\xi(t)= g^{-1}\dot g(t) \in
se(n)$.  In this case $\omega=R^{-1}\dot R(t)$ and $v=R^{-1}\dot x(t)$
are respectively the angular velocity matrix and the vector of linear
velocity of $A$ in the body frame.

As in the classical case, we suppose that the center of mass ${\bf C}$
of the body does not coincide with the origin $A$ of the body frame.
Let $(a_1,\dots, a_n)^T$ be constant position vector of ${\bf C}$ in
this frame and, as above, $J=\mbox{diag} (J_1, \dots, J_n)$ be its
mass tensor.  Then the Lagrangian is
\begin{gather}
L= - \frac 14 \mbox{tr } (\xi \, {\mathbb J}\, \xi^T )
\equiv - \frac 14 \mbox{tr } (\omega (J+m a\otimes a) \omega)+ m (v, \omega a) + \frac m2 (v,v) ,
\label{cont.sl.lagr} \\
{\mathbb J} = S \; \textup{diag } (J_1, \dots, J_n, m)\; S^T, \quad 
S=\begin{pmatrix} 1 &  &  & a_1 \\
                     & \ddots & & \vdots \\
                     &        & 1 & a_n \\
                   0 & \dots & 0 & 1   \end{pmatrix}\in SE(n), \nonumber
\end{gather}
where $S$ describes the position of the center of mass ${\bf C}$
relative to the body frame. 


The body momentum is an element of the dual space $se ^{*}(n)$ and it is given by the pair
\begin{gather*}
P=(M, p)\in se^*(n), \qquad M\in so(n), \quad p\in {\mathbb R}^n, \\
M_{ij} = \pder{L}{\omega_{ij}}, \quad p_i=\pder {L}{v_i}, \qquad i,j=1,\dots,n.
\end{gather*} 
Straightforward evaluation leads to the formulae
\begin{equation} 
\label{cont.mom.se(n)}
\begin{aligned}
M & =(J+m a\otimes a) \omega +\omega (J+m a\otimes a)+ m(v\otimes a-a\otimes v)\in so(n), \\
p  &= m (v+ \omega a) \in{\mathbb R}^n .
\end{aligned}
\end{equation} 
Here $M$ is the angular momentum of the body with respect to its center of mass $C$ and
$p$ is the linear momentum of $C$ as a point with mass $m$. 

\paragraph{Left-invariant constraints on $SE(n)$.} 
There are numerous ways to introduce nonholonomic constraints for the
generalized Chaplygin sleigh. For example, one can require that the
velocity of the reference point is restricted to a $k$-dimensional
linear subspace fixed in the body. For $n=3$, such constraints were
studied in \cite{NeFu} and \cite{ZBl2002}.)

Another natural choice is to define the constraint subspace
${\mathfrak d}\in se(n)$ to be the set of matrices of the form
\begin{equation} \label{d_for_se_n}
\begin{pmatrix}
0 &  \omega_{12} & \cdots & \omega_{1n} & v_1 \\ 
-\omega_{12} &  0 &        & 0           & 0 \\ 
\vdots       &   &  &   & \vdots  \\
- \omega_{1n} & 0  &  & 0  & 0 \\
0 & 0 & \dots & 0 & 0 
\end{pmatrix} . 
\end{equation}
In the particular case $n=2$ we have
\begin{equation} \label{se2.constr}
S=\begin{pmatrix}
0 & -\omega & v_1 \\ 
\omega & 0 & 0 \\ 
0 & 0 & 0 \end{pmatrix} .
\end{equation}

\section{Discrete Suslov System on $SO(n)$}
Now we apply the  the discrete Euler--Poincare-- Suslov equations 
(\ref{disc.momentum.eqn}) 
to construct a discretization of the Suslov problem. 
Let $R_k\in SO(n)$ be the orthogonal rotation matrix describing the $k$-th position 
of $n$-dimensional top. 

Introduce the {\it finite rotation matrix\/} 
$\Omega_k = R_{k}^T R_{k+1}$, analog of the angular velocity $\omega$ in the body. 
Note that in the continuous limit, when $R_{k+1}= R_{k} + \varepsilon \dot R$, 
$\varepsilon <<1$, one has
\begin{equation} \label{limits}
\Omega_k = {\bf I} + R^{-1} \dot R = {\bf I}+ \varepsilon \omega,
\end{equation}

Define the left-invariant discrete Lagrangian on $SO(n)\times SO(n)$ by substituting $\omega$
in (\ref{cont.so-n.lagr}) by $ R_{k}^T ( R_{k+1}-R_{k})\equiv \Omega_k-{\bf I}$. 
Using the property $R_{k}^T R_{k}={\bf I}$, we get
$$
l_d (\Omega_k) = \frac 12 \operatorname{tr} (\Omega_k J), \quad
 \mbox{and} \quad L_d (R_{k}, R_{k+1})=\frac 12 \operatorname{tr} (R_k J R_{k+1}^T) .
$$ 
Then, following the definition (\ref{disc.mom.dfn}), the body angular momentum 
$M_k \in so^*(n)$ has the form
\begin{equation} \label{*}
M_k = R_k^{-1} R_{k+1} J - J R_{k+1}^T R_k \equiv \Omega_k J - J\Omega_k^T ,
\end{equation}
which in the above limit transforms to $J\omega +\omega J$, the standard relation between
the angular velocity and momentum. 
The expressions for $L_d, M_k$ were originally introduced in \cite{MosVes}. 

\paragraph{Remark.} In the classical case $n=3$ one can parameterize $R_k$
in terms of the Euler angles $\theta_k,\psi_k,\phi_k$, as coordinates on $SO(3)$
(see e.g., \cite{Whitt}),
$$
\begin{pmatrix}
\cos{\phi }_k \cos {\psi }_k - \cos {\theta }_k \sin {\phi }_k \sin {\psi }_k 
& - \cos {\phi }_k \sin {\psi }_k - \cos {\theta }_k \sin {\phi }_k \cos {\psi }_k
& \sin {\theta }_k \sin {\phi }_k \\   
\sin {\phi }_k \cos {\psi }_k +   \cos {\theta }_k \cos {\phi }_k \sin {\psi }_k 
&    - \sin {\phi }_k \sin {\psi }_k + \cos {\theta }_k \cos {\phi }_k \cos {\psi }_k 
&    - \sin {\theta }_k \cos {\phi }_k \\   
\sin {\theta }_k \sin {\psi }_k & \sin {\theta }_k  \cos {\psi }_k& \cos {\theta }_k
\end{pmatrix}
$$
Substituting these ones and analogous expressions for $R_{k+1}$ into the discrete Lagrangian 
$L_d (R_k, R_{k+1})$, we obtain
\begin{align} 
 L_d & = \frac 12 \bigg[ \cos \theta_k\,\cos \theta_{k + 1}\,
 [ 1 + \cos (\Delta\phi_k)\,\cos (\psi_k + \psi_{k + 1})] \nonumber \\
  & \quad + \cos (\psi_k + \psi_{k + 1})\,\sin \theta_k\,\sin{\theta}_{k + 1} 
+\cos(\Delta\phi_k) [ \sin \theta_k \,\sin \theta_{k + 1} \nonumber \\
& \quad -\cos(\psi_k +\psi_{k + 1}) ]
 +\frac 12 [\cos \theta_{k+1}-\cos {\theta}_{k}] \,\sin (\Delta{{\phi }_k})\,
      \sin ({{\psi }_k} + {{\psi }_{k + 1}})\bigg] A_1 \nonumber \\ 
& \quad + \frac 12 \bigg[\cos \vartheta_k\,\cos\vartheta_{k + 1} +
\cos (\Delta \phi_k)\,\cos(\Delta\psi_k)  
 - \cos \vartheta_k\,\cos \vartheta_{k + 1}\,\cos(\Delta\phi_k)\,
\cos (\Delta\psi_k) \nonumber \\
& \quad + \cos (\Delta\phi_k)\,\sin{{\vartheta}_k}\,\sin\vartheta_{k + 1}
 -\cos (\Delta\psi_k)\,\sin\vartheta_k\,\sin\vartheta_{k + 1} \nonumber \\ 
 & - \cos\vartheta_k\,\sin ({{\phi }_k} + {{\phi }_{k + 1}})\,\sin(\Delta\psi_k)
- \cos \vartheta_{k + 1}\,\sin(\Delta\phi_k)\,
\sin({{\psi }_k} + {{\psi }_{k + 1}})\bigg] A_2 \nonumber \\
& \quad - \frac 12 \bigg[  \cos(\Delta\phi_k)\,\cos(\Delta\psi_k) +
\cos \theta_k\,\cos \theta_{k + 1}\,\cos(\Delta\phi_k)\,\cos (\Delta\psi_k) 
\nonumber \\ 
  &\quad  - \cos \theta_k\,\cos\theta_{k + 1}  
+ \sin\theta_k\,\sin\theta_{k + 1}( \cos(\Delta\psi_k)- \cos(\Delta\phi_k)) 
\nonumber \\
& \quad - (\cos \theta_k +\cos {\theta}_{k + 1})
\sin(\Delta\phi_k)\,\sin (\Delta\psi_k)\bigg] A_3 \label{L_angles}
\end{align}
where $\Delta\theta_k =\theta_{k + 1}-\theta_k$, 
$\Delta\phi_k =\phi_{k + 1}-\phi_k$, $\Delta\psi_k =\psi_{k + 1}-\psi_k$, and
$$
A_1=J_2+J_3, \quad A_2=J_1+J_3, \quad A_3=J_1+J_2
$$
are the principal moments of inertia of the rigid body.  

In the continuous limit, setting in (\ref{L_angles})
\begin{equation} \label{angle_limit}
\theta_{k+1}-\theta_k =\dot\theta\, \delta t, \quad \phi_{k+1}-\phi_k =\dot\phi\,\delta t, 
\quad \psi_{k+1}-\psi_k =\dot\psi\, \delta t, \qquad \delta t << 1
\end{equation}
then expanding in $\delta t$ and dividing by $(\delta t)^2$, up to an additive constant 
and terms of order $\delta t$, one obtains 
the well-known expression for the kinetic energy of the top (see, e.g., \cite{Whitt}) 
\begin{align} 
T & = \frac 12 (\dot\phi\sin\theta\sin\psi+\dot\theta\cos\psi)^2A_1 \nonumber \\ 
 & \qquad +\frac 12 (\dot\phi\sin\theta\cos\psi -\dot\theta\sin\psi)^2 A_2 + 
\frac 12 (\dot\psi + \dot\phi\cos \theta)^2 A_3,  \label{lagr_cont}
\end{align}
where the expressions in brackets represent components of the angular velocity vector
in the frame attached to the body.

Notice that the discrete Lagrangian (\ref{L_angles}) does not coincide with 
the "straightforward" discretization of (\ref{lagr_cont}) obtained with a direct replacement of  
the angular velocities by the angular differences according to (\ref{angle_limit}). 
\medskip

\paragraph{Discrete constraints on $SO(n)$.} Following the approach described in 
Section 2, we impose discrete left-invariant constraints on $SO(n)\times SO(n)$
in the form of restrictions on finite rotations $\Omega_k \in SO(n)$.
In accordance with the continuous constraints (\ref{Suslov_constraints}),  
we assume that admissible rotations must be exponents of the vectors of the linear space
$$
{\mathfrak d}= \mbox{span} \{ e_{1} \wedge e_n,\; \dots, e_{n-1}\wedge e_n \}
\subset so(n) . 
$$

\begin{lemma} \label{symm} \begin{description}
\item{1).} In the basis $e_1,\dots, e_n$,
the admissible rotation matrices have the structure
\begin{equation} \label{sp}
(\Omega_k )_{ij} = (\Omega_k)_{ji}, \quad 
(\Omega_k )_{in} = -(\Omega_k)_{ni}, \qquad 1 \le i, j\le n-1,
\end{equation}
that is, they are anti-symmetric in its last row and column 
and symmetric in the other part. 

\item{2).} The admissible displacement subvariety $S=\exp{\mathfrak d}$ is 
homeomorphic to the projective space ${\mathbb P}^{n-1}=S^{n-1}/{\mathbb Z}^2$. 
In the same basis, the components of $\Omega_k$ are parameterized by points 
of the unit sphere $S^{n-1}=\{z_0^2+z_1^2+\cdots +z_{n-1}^2=1\}$ in the form
\footnote{Here and below, to simplify notation, we omit the discrete time index $k$ at 
the components of $z$.}

\begin{gather} \label{sp1}
(\Omega_k )_{ij}=\delta_{ij} -2 z_i z_j, \quad
(\Omega_k )_{in} = -(\Omega_k)_{ni} = 2z_0 z_i, \quad 
(\Omega_k )_{nn}= 2z_0^2-1 ,\\
1 \le i, j\le n-1 . \nonumber
\end{gather}
\end{description}
\end{lemma}

Note that in the continuous limit described by (\ref{limits}), 
conditions (\ref{sp}) yield the constraints (\ref{matrix_constraint}) on $so(n)$. 
\medskip

\noindent{\it Proof of Lemma } \ref{symm}. 
1). Any vector of ${\mathfrak d}\subset so(n)$ can be represented in the form
$\theta {\bf u}\wedge e_n$, where $\theta$ is a nonzero constant and 
${\bf u}=(u_1,\dots, u_{n-1},0)^T$ is a
unit vector in ${\mathbb R}^{n-1}=\mbox{span }(e_1,\dots, e_{n-1})$.
The odd powers of $\theta {\bf u}\wedge e_n$ are skew-symmetric and 
have zero left-upper $(n-1)\times (n-1)$ part, whereas the even powers 
are symmetric and have zero last row and last column. Hence, the exponent of
$\theta {\bf u}\wedge e_n$ must be of the form (\ref{sp}).
 
2). The operator ${\cal R}_{\theta, u}=\exp (\theta {\bf u}\wedge e_n)\subset SO(n)$ 
describes rotation in the 2-plane spanned by ${\bf u}, e_n$ by the angle $\theta$. 
Then we get
\begin{align*}
{\cal R}_{\theta, u} \,{\bf u} & = \cos\theta\cdot{\bf u}- \sin\theta \cdot e_n, \\ 
{\cal R}_{\theta, u} \, e_j & = e_j- (e_j, {\bf u} ){\bf u} 
+ (e_j, {\bf u}) (\cos\theta\cdot{\bf u}- \sin\theta \cdot e_n),
\qquad 1 \le j\le n-1 , \\
{\cal R}_{\theta, u}\, e_n & = \cos\theta \cdot e_n +\sin\theta\cdot {\bf u} \, .
\end{align*}
The latter $n$ vectors form columns of the matrix ${\cal R}_{\theta, u}$. 
Setting in the above formulas
\begin{equation} \label{t->q}
z_i= \sin\theta/2 \, u_i, \quad 1 \le i, j\le n-1 , \quad z_0=\cos\theta/2 
\end{equation}
and identifying ${\cal R}_{\theta, u}$ with $\Omega_k$ we arrive at expressions (\ref{sp1}). 

Since $\sin(\theta/2)= - \sin(2\pi-\theta)/2$ and $\cos(\theta/2)= - \cos(2\pi-\theta)/2$, 
from (\ref{t->q}) we conclude that opposite points on $S^{n-1}$ correspond to the
same admissible rotation ${\cal R}_{\theta, u}$. 
Finally, there is a bijection between $S=\exp{\mathfrak d}$ and 
${\mathbb P}^{n-1}=S^{n-1}/{\mathbb Z}^2$. 
The lemma is proved. 
\medskip

Note that (\ref{sp}) imply left-invariant constraints on $SO(n)\times SO(n)$ in the form
$$
\operatorname{tr} (R_k^T \, e_j\wedge e_n\, R_{k+1}-R_{k+1}^T\, e_j\wedge e_n\, R_k)=0\, , 
\qquad j=1,\dots,n-1. 
$$
\medskip

\paragraph{Rotations about an axis.} In the classical case $n=3$  
the conditions (\ref{sp}) say that $\Omega_k$ is a finite rotation about an axis 
lying in the plane $(e_1, e_2)$, while expressions (\ref{sp1}) imply that the rotation
axis is directed along vector $\rho=(z_2, -z_1, 0)^T\in {\mathbb R}^3$.

Indeed, the group $SO(3)$ can be regarded as covered twice by 
the unit sphere $S^3=\{q_0^2+q_1^2+q_2^2+q_3^2=1\}$, where $q_0,\dots,q_3$ are
the Euler--Rodriguez parameters such that
any rotation matrix $W\in SO(3)$ can be represented in form (see, e.g., \cite{Whitt})
\begin{equation} \label{Rod}
W = \begin{pmatrix} 
q_{0}^{2}+q_{1}^{2}-q_{2}^{2}-q_{3}^{2} & 2(q_{1}q_{2}+q_{3}q_{0}) & 
-2(q_{1}q_{3}-q_{2}q_{0}) \\ 
2(q_{1}q_{2}-q_{3}q_{0}) & q_{0}^{2}+q_{2}^{2}-q_{1}^{2}-q_{3}^{2} & 
-2(q_{2}q_{3}+q_{0}q_{1}) \\ 
-2(q_{1}q_{3}+q_{2}q_{0}) & -2(q_{2}q_{3}-q_{0}q_{1}) & 
q_{0}^{2}+q_{3}^{2}-q_{1}^{2}-q_{2}^{2}
\end{pmatrix} .
\end{equation}
The operator $W$ describes a finite rotation in ${\mathbb R}^3$ 
about the vector 
${\bf e} = (q_1, q_2, q_3)^T$ by the angle $\theta$ such that $q_0=\cos \theta/2$.

Setting in (\ref{Rod}) $(\Omega_k)_{12}=(\Omega_k)_{21}$
implies $q_3=0$, hence $W$ is a rotation about an axis lying in the plane 
$(e_{1}, e_{2})$. In this case admissible operators $\Omega\in {\cal S}\subset SO(3)$ 
have the form 
\begin{equation} \label{PE}
\Omega= \begin{pmatrix} 2(q_0^2+q_1^2)-1 & 2q_1 q_2 & 2q_0 q_2 \\
2 q_1 q_2 & 2(q_0^2+q_2^2)-1 & -2q_0 q_1 \\
- 2q_0 q_2 & 2 q_0 q_1 & 2q_0^2 -1 \end{pmatrix},
\end{equation}
which, under the substitution $q_1=- z_2$, $q_2=z_1$, $q_0=z_0$, coincides with
the parameterization (\ref{sp1}).
As a result, the variety of such matrices is the real projective plane 
$\mathbb{RP}^2 =S^2/{\mathbb Z}^2$. 

We emphasize that, in general, the $k$-th position of the body
$R_k=\Omega_{k-1} \cdots \Omega_0$
{\it is not \/} a rotation in the plane $(e_{1}, e_{2})$.
\medskip

\paragraph{Discrete momentum locus ${\cal U}\subset so^*(3)$.}
In contrast to the continuous case, 
the discrete momentum $M_k$ does not lie in a linear subspace in the coalgebra $so^*(3)$,
but on a nonlinear algebraic variety ${\cal U}\subset so^*(3)$ defined by 
the relation (\ref{*}) and the conditions (\ref{PE}). 

If in the frame $e_1, e_2, e_3$ the tensor $J$ is diagonal,  
$J=\operatorname{diag}(J_1, J_2, J_3)^T$, then
the angular momentum vector 
$ \overrightarrow {M} =(M_1= -M_{23}, M_2= M_{13}, M_3= -M_{12})^T$  has the form 
$$
\overrightarrow {M}= 2( (J_2+J_3)q_0 q_1,(J_1+J_3)q_0 q_2, (J_1 - J_2)q_1 q_2)^T
$$
(as above, to avoid tedious notation we omit the discrete time index at the components 
of $q$). Here and below, without loss of generality, we always assume $q_0 \ge 0$. 
As a result, ${\cal U}$ coincides with {\it the Steiner Roman surface} in ${\mathbb R}^3$ 
given by the quartic equation
\begin{gather}
\frac{J_1- J_2}{(J_2+J_3)(J_1+J_3)}  M_1^2 M_2^2 
+\frac{J_1+J_3}{(J_2+J_3)(J_1-J_2)}  M_1^2  M_3^2 \nonumber \\
+\frac{J_2+J_3}{(J_1+J_3)(J_1-J_2)}  M_2^2 M_3^2 
- 2 M_1  M_2  M_3=0 \label{Steiner}
\end{gather}
(see, e.g., \cite{Gray, No}). 

In general case, when $J$ is not diagonal in this frame, one has the parameterization
\begin{equation} \label{full}
\overrightarrow {M} =
2 \begin{pmatrix} (J_{22}+J_{33}) q_0 q_1 - J_{12}q_0q_2 & - (J_{13}q_1 +J_{23}q_2)q_2 \\
(J_{11}+J_{33}) q_0 q_2  - J_{12} q_0 q_1 & + (J_{13} q_1  + J_{23} q_2) q_1 \\
(J_{11}-J_{22}) q_1 q_2  - J_{12} (q_1^2-q_2^2) & - (J_{13}q_1 + J_{23}q_2) q_0 
\end{pmatrix} .
\end{equation}
One can show that the components of $\overrightarrow M$ satisfy an
algebraic equation of degree 4, which generalizes (\ref{Steiner}) and
which we do not write here.  The corresponding algebraic surface $\cal
U$ in ${\mathbb R}^3=(M_1,M_2, M_3)$ has pinch points and
self-intersections. One can also show that if the quadratic form
$(J_{22}+J_{33})q_{1}^{2} -2J_{12} q_{1}q_{2} +( J_{33}+J_{11})
q_{2}^{2}$ is positive-definite, then any pair $(M_1,M_2)$ has at most
two real inverse images on $\cal U$.

\begin{figure}[h,t]
\begin{center}
\includegraphics[width=.6\textwidth]{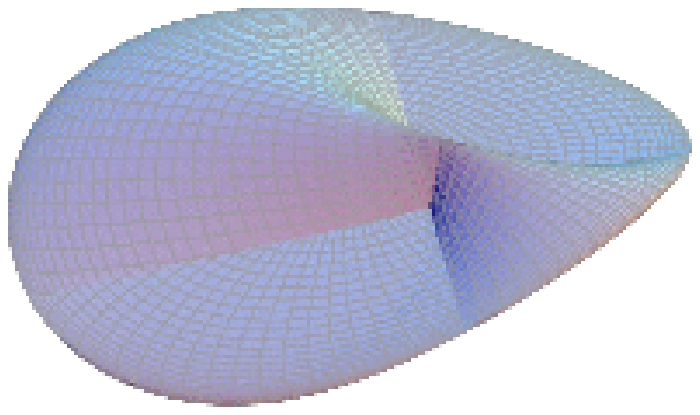}
\caption{\footnotesize The Momentum Surface $\cal U$.}
\label{st.fig}
\end{center}
\end{figure}

\begin{figure}[h,t]
\begin{center}
\includegraphics[width=.6\textwidth]{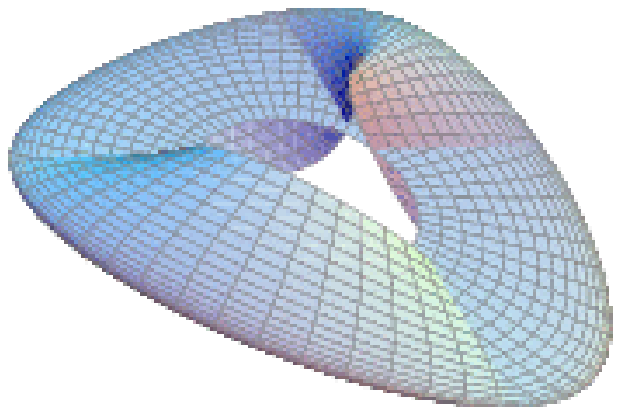}
\caption{\footnotesize A Circular Section of the Momentum Surface}
\label{st1.fig}
\end{center}
\end{figure}

An example of such a surface for an unbalanced
inertia tensor and its circular section for $0.4\le q_0 \le 0.8$ are given in Figures 
\ref{st.fig}, \ref{st1.fig} respectively. \medskip

\medskip

\paragraph{Discrete EPS equations on $so^*(3)$.} In the considered case $G=SO(3)$,
the discrete momentum equation with multipliers (\ref{disc.momentum.eqn}) takes the form
\begin{equation} \label{suslov}
M_{k+1} = \Omega_{k}^T M_{k} \Omega_{k} + \lambda_{k} \,
\begin{pmatrix}  0 & 1 & 0 \\
                 -1 & 0 &0 \\
                   0 & 0 & 0 \end{pmatrix}, \quad
M_k = \Omega_k J - J\Omega_k^T ,
\end{equation}
where the components of $\Omega_k$ are subject to constraints (\ref{sp}). 

This provides a discrete analog of the Suslov system (\ref{Suslov_eq}) on $so^*(3)$ 
and defines a map ${\cal U}\to {\cal U}$ or, in view of expressions (\ref{PE}), 
(\ref{full}), a map \\ 
${\cal B}\: :\, {\mathbb RP}^2 \to {\mathbb RP}^2\; :\;
(q_1,q_2,q_0) \to (\tilde q_1,\tilde q_2, \tilde q_0)$, which is generally
multi-valued. 
\medskip

To describe the latter map in details, we note that in (\ref{suslov})
\begin{align} \overrightarrow { \Omega_{k}^T M_{k} \Omega_{k}}& 
\equiv \overrightarrow{J\Omega_k -\Omega_k^T J}\equiv \Omega_{k}^T  \overrightarrow {M_k}
\nonumber \\
& = 2 \begin{pmatrix} (J_{22}+J_{33}) q_0 q_1  - J_{12}q_0q_2 & + (J_{13}q_1 + J_{23}q_2)q_2 \\
(J_{11}+J_{33}) q_0 q_2  - J_{12} q_0 q_1  & - (J_{13} q_1 + J_{23} q_2) q_1 \\
-(J_{11}-J_{22}) q_1 q_2  + J_{12} (q_1^2-q_2^2) & - (J_{13}q_1 + J_{23}q_2) q_0 
\end{pmatrix} , \label{full-1}
\end{align}
where $\overrightarrow\omega$ denotes vector representation of element $\omega$ of $so(3)$.
Comparing this with (\ref{full}), we find that (\ref{suslov}) can be written in form
\begin{equation} 
\label{diff_suslov}
\overrightarrow {M_{k+1}} = \overrightarrow {M_{k}} +
4 \begin{pmatrix}  (J_{13}q_1 + J_{23}q_2)q_2 \\
 - (J_{13} q_1 + J_{23} q_2) q_1 \\
-(J_{11}-J_{22}) q_1 q_2  + J_{12} (q_1^2-q_2^2) + \lambda_k \end{pmatrix} ,
\end{equation}
which can be viewed as a discrete analog of equations (\ref{ep3.26}). This also shows
that the difference vector $\overrightarrow {M_{k+1}} -\overrightarrow {M_{k}}$
is orthogonal to the rotation axis directed along $(q_1,q_2,0)\in {\mathbb R}^3$,  
as expected. 

As a result, the map ${\cal B}\: :\, \mathbb{RP}^2 \to \mathbb{RP}^2$ given by
(\ref{suslov}) consists of the following 3 steps:

\begin{description} 
\item{1).} Given original set $q_1, q_2, q_0= \sqrt{1-q_1^2-q_2^2}$, one finds components 
of $\overrightarrow {M_k}$ from (\ref{full}) and of 
$\overrightarrow {\Omega_{k}^T M_{k} \Omega_{k}}$ from (\ref{full-1}).

\item {2).} Given the components  
$$
(\overrightarrow {M_{k+1}})_1 =(\overrightarrow {\Omega_{k}^T M_{k} \Omega_{k}})_1, \quad 
(\overrightarrow {M_{k+1}})_2 =(\overrightarrow {\Omega_{k}^T M_{k} \Omega_{k}})_2,
$$
one finds new $\tilde q_1, \tilde q_2$ by 
solving the system of two algebraic equations originating from (\ref{full}) 
\begin{equation} 
\begin{aligned} 
( \overrightarrow {M_{k+1}})_1 & = \left ( (J_{22}+J_{33})\tilde q_1 
- J_{12}\tilde q_2 \right) \sqrt{1-\tilde q_1^2-\tilde q_2^2}
 - (J_{13}\tilde q_1 +J_{23}\tilde q_2)\tilde q_2,  \\
(\overrightarrow {M_{k+1}})_2 & = \left( (J_{11}+J_{33})\tilde q_2
-J_{12}\tilde q_1\right ) \sqrt{1-\tilde q_1^2-\tilde q_2^2}
 + (J_{13}\tilde  q_1  + J_{23}\tilde q_2)\tilde q_1 .
\end{aligned} \label{sys}
\end{equation}
In ${\mathbb R}^3=(q_1, q_2, q_0)$ these equations describe two centrally symmetric 
quadratic surfaces $Q_1, Q_2$ which intersect
the unit sphere $q_1^2+q_2^2+q_0^2=1$ along curves $C_1, C_2$ respectively. 
Each curve is a union of two  ovals, which are centrally 
symmetric to each other. 
The intersection of $C_1, C_2$ gives 4 complex points and 
2 or none real points on ${\mathbb P}^2$.
Thus there are at most two different real solutions 
$(\tilde q_1^{(j)} , \tilde q_2^{(j)} ,\tilde q_0^{(j)})$ with $\tilde q_0^{(j)}>0$.

\item{3).} One chooses a solution $\left(\tilde q_1^{(1)} , \tilde q_2^{(1)}
,\tilde q_0^{(1)}>0\right)$ and finally finds the last component
$(\bar M_{k+1})_3$ by the formula
$$ (\bar M_{k+1})_3=(J_{11}-J_{22}) \tilde q_1 \tilde q_2 - J_{12}
(\tilde q_1^2-\tilde q_2^2)  - (J_{13}\tilde q_1 + J_{23}\tilde q_2)
\tilde q_0,
$$ which is obtained from (\ref{sys}) by substitutions $k \to k+1$ and
the $q \to \tilde q$.
\end{description}
\medskip

As a result, for $n=3$ the map $M_k\to M_{k+1}$ given by (\ref{suslov}) 
is generally {\it 4-complex valued} and 2-real valued. 
In order to choose one of the 2 real branches, 
we must use some extra arguments, like existence of an additional integral,
or, at least, to restrict ourselves with sufficiently small $q_1, q_2$, which correspond to 
rotations $\Omega$ by a small angle $\theta$. In this case only one of the solutions
$(\tilde q_1^{(j)} , \tilde q_2^{(j)})$ will be small and it is natural to choose it.
\medskip

It appears that the constrained energy integral (\ref{sus_int_1}) of the continuous
Suslov system is preserved by the discrete system as well. 

\begin{theorem} \label{ints}
The discrete Suslov system (\ref{suslov}) has quadratic integral 
\begin{equation} \label{sus_int}
( J_{11}+J_{33}) M_1^2 +2 J_{12} M_1 M_2+ ( J_{22}+J_{33}) M_2^2 ,
\end{equation} 
which gives rise to the following quartic integral
 in terms of the parameters $q_0,q_1,q_2$:
\begin{align} 
H= ( ( J_{22}+J_{33})q_{1}^{2} & -2J_{12} q_{1}q_{2} 
+( J_{11}+J_{33}) q_{2}^{2} ) \nonumber \\ 
& \cdot
\left( ( J_{13}q_{1} +J_{23} q_{2})^{2}+ [(J_{11}+J_{33})(J_{22}+J_{33})-J_{12}^{2} ]
q_0^{2} \right) . \label{quatern_int}
\end{align} 
\end{theorem}

The proof is straightforward: substituting expressions (\ref{full}) and (\ref{full-1}) into
(\ref{sus_int}) gives the same expression in terms of $q_0,q_1,q_2$.
\medskip 

The fact that (\ref{sus_int}) does not depend on $M_3$ is quite natural:
different branches of the map (\ref{suslov}) have the same value of the integral.

It should be  emphasized that the complete energy integral $(M, {\mathbb I}^{-1} M)$ of 
the continuous Suslov problem is not preserved in the discrete setting.

\paragraph{Invariant curves.}
As follows from Theorem \ref{ints},  the map has invariant curves, 
which are either intersections of the sphere 
$\{q_1^2+q_2^2 + q_0^2=1\}$ with a quartic surface $H(q)=h$ or, in the momentum space $so^*(3)$, 
intersections of the  generalized quartic Steiner surface $\cal U$ with 
elliptic cylinders defined by (\ref{sus_int}).
Thus,  the invariant varieties are algebraic curves of order 8. 

Assume that quadratic form $(J_{22}+J_{33})q_{1}^{2} -2J_{12} q_{1}q_{2} 
+( J_{33}+J_{11}) q_{2}^{2}$ is positive definite. Then, as follows from 
(\ref{quatern_int}),  on the upper hemisphere $0\le q_0\le 1$ 
real invariant curves 
consist of two branches: for small positive values of $h$ one branch is a small oval around the
origin $(0,0)$ whereas the other branch is an oval close to the equator $\{q_0=0\}$ of the sphere.
It may or may not intersect  the equator. In the first case the opposite points of 
intersection are identified.

These different branches correspond to the two connected components of the intersection
of the Steiner surface $\cal U$ with the cylinder. 

\begin{figure}[h,t]
\begin{center}
\includegraphics[width=0.5\textwidth]{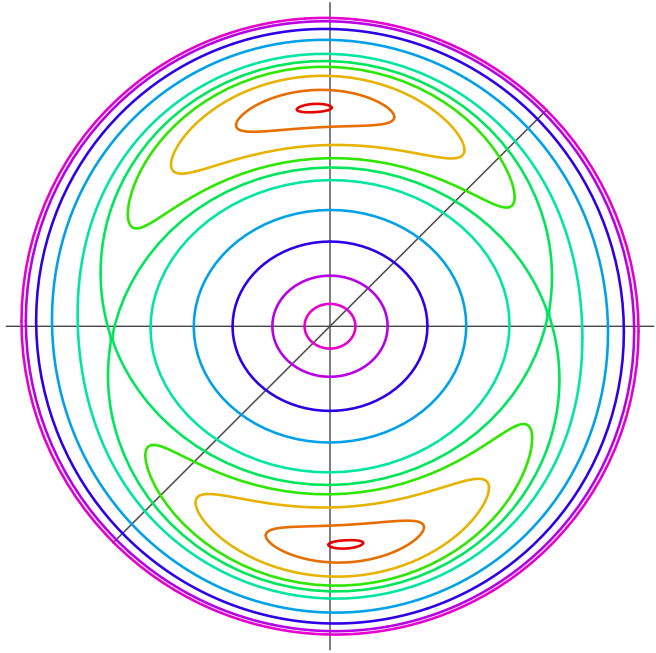}
\caption{\footnotesize Invariant Curves and the Equilibria Line on $\mathbb{RP}^2$.}
\label{quat.fig}
\end{center}
\end{figure}

As value of the integral increases, the branches approach each other: the
smaller one becomes bigger
and the bigger shrinks. At a certain critical value $h=h^*$ the branches intersect at two
opposite saddle points and form a separatrix, and for the next critical value $h^{**}>h^*$ 
the two branches shrink to opposite center points. There are no real invariant curves
for $h>h^{**}$. Note that for $h=h^*$ and $h=h^{**}$ 
the elliptic cylinder is tangent to the surface $\cal U$. 
An example of the invariant curves foliation is given in Figure \ref{quat.fig}. 
\medskip

\noindent{\bf Remark.} As noticed in \cite{MosVes}, in the absence of nonholonomic constraints, 
the map $M_k \to M_{k+1}$ given by the
discrete Euler--Poincar\'e equations (\ref{suslov}) is multi-valued, because, in general,
the equation $M_k=\Omega_k J - J\Omega_k^T$ has more than one solution. 

In presence of the constraints (\ref{PE}), 
the latter equation has generally {\it a unique} solution 
(except the points on self-intersection on $\cal U$), 
however, as we saw above, the choice of $\lambda_{k+1}$ or $(\bar M_{k+1})_3$ is not
unique, and the map describing the discrete Suslov problem is multi-valued as well. 
\medskip

\paragraph{Stationary solutions of the discrete Suslov problem.} 
As follows from (\ref{diff_suslov}), 
if the initial values $q_1, q_2$ satisfy the condition $J_{13} q_1 + J_{23} q_2=0$, then
\begin{gather*}
(\overrightarrow { M_{k+1}})_1 =(\overrightarrow { M_{k}})_1, \quad 
(\overrightarrow {M_{k+1}})_2 =(\overrightarrow {M_{k}})_2, \quad
(\overrightarrow  {\Omega_{k}^T M_{k} \Omega_{k}})_3 = - (\overrightarrow { M_{k}})_3,
\end{gather*}
that is, 
the coadjoint action $M_k\mapsto \Omega_{k}^T M_{k} \Omega_{k}$ is the mirror reflection
with respect to the plane $M_3=0$. Then it is natural to choose the multiplier $\lambda_k$
such that $ (\bar M_{k+1})_3=(\bar M_{k+1})_3$. 

As a result, {\it one of the branches} of the map ${\cal B}$ 
has a one-parametric family of stationary solutions (equilibria) 
characterized by points of the line
$$
{\mathbb P}= \{ S^2\cap \{J_{13} q_1 + J_{23} q_2=0 \} \}/{\mathbb Z}^2.
$$ 
They correspond to discrete versions of permanent rotations of the body in the classical
Suslov problem. (In Figure \ref{quat.fig} the set of equilibria points
is represented by a straight line segment.) 

In view of (\ref{PE}), opposite points $(q_1, q_2, q_0)$ and $(-q_1, -q_2, q_0)$ 
on ${\mathbb P}$ correspond to mutually inverse finite rotations 
$\Omega$ and $\Omega^T$ respectively. 

As also follows from (\ref{diff_suslov}), there are no equilibria points outside of this line.
In particular, neither the saddle points nor the centers of the invariant foliation
on ${\mathbb RP}^2$ are stationary points. 

Finally, note that, like in the continuous system, for a balanced inertia tensor $J_{13}=J_{23}=0$
{\it all} the solutions of (\ref{diff_suslov}) are stationary, i.e., the discrete body momentum
$M_k$ is preserved. 
\medskip

\noindent{\bf Remark.}  The foliation of $\mathbb{RP}^2$ by invariant
curves gives us a natural way of choosing the branches of the map $\cal B$
in the general case.  Namely, if the initial point $(q_1, q_2)$ lies
in the domain $\cal S\subset {\mathbb{RP}}^2$ defined by the
condition $0< h\le h^*$, i.e., it represents either a relatively small
or sufficiently big finite rotation $\Omega$, then the points $(q_1,
q_2)$ and $(\tilde q_1, \tilde q_2)$ have to belong to the same connected
component of the invariant curve.  In other words, if the initial
point lies in the interior (exterior) part of $\cal S$, one has to
choose a real solution of (\ref{sys}) that has the smallest
(largest) norm $\tilde q_1^2 +\tilde q_2^2$, respectively.

On the other hand, if $(q_1, q_2)$ lies in complement
$\mathbb{RP}^2\setminus {\cal S}$, i.e., it is between the
separatrices, then a real initial point $(q_1, q_2)$ may lead to
complex $(\tilde q_1, \tilde q_2)$ only. In particular, when the
initial point is a center, the next point is necessarily complex,
although the value of the integral remains to be real.

If branches of the map $\mathbb{RP}^2\mapsto \mathbb{RP}^2$ are
chosen according to the above way, then the discrete time dynamics
inherits all the main properties of the continuous Suslov problem.

Namely, let $\Delta_{-}$ and $\Delta_{+}$ denote semi-planes of
$\mathbb{RP}^2$ defined by conditions $J_{13} q_1 + J_{23} q_2 <0$
(respectively $>0$) and let $\Theta_{-}$ and $\Theta_{+}$ be
semi-planes given by
$$
(J_{12}J_{13}+J_{22}J_{23} +J_{23}J_{33}) q_{1}
-( J_{11}J_{13}+J_{12}J_{23}+J_{13}J_{33}) q_{2} < 0 ,
$$
respectively $>0$.

\begin{theorem} \label{discr_sus_dyn}
If the initial point ${\bf q}=(q_1, q_2)$ lies in the interior part of
${\cal S}\subset \mathbb{RP}^2$, then for $k \to -\infty$ and $k \to
+\infty$
the sequence $\{{\bf q}_k\}$ remains on the same branch of invariant curve
and tends to the unstable equilibria semi-line
${\mathbb P}_u ={\mathbb P} \cap \Theta_{-}$ and the
stable equilibria semi-line ${\mathbb P}_s={\mathbb P} \cap \Theta_{-}$
respectively. It lies entirely in one of the semi-planes $\Delta_{\pm}$.
\end{theorem}

\begin{figure}[h,t]
\begin{center}
\includegraphics[width=0.5\textwidth]{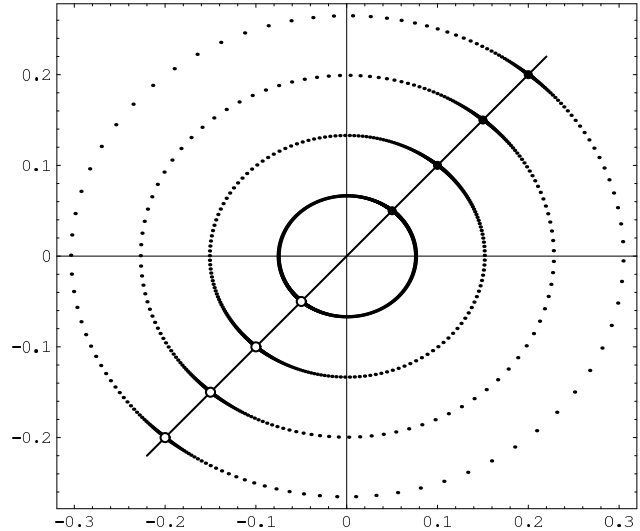}
\caption{\footnotesize Discrete Dynamics near the origin of $\mathbb
{RP}^2$.}
\label{quat_dots}
\end{center}
\end{figure}

For the foliation indicated in Figure \ref{quat.fig}, the corresponding
discrete time dynamics
in the neighborhood of the origin is
given in Figure \ref{quat_dots}, where stable and unstable equilibria points
on
$\mathbb{P}$ as depicted as dots and circles respectively.

\vspace{2em}
As follows from Theorem \ref{discr_sus_dyn}, for $k \to -\infty$ and $k \to
+\infty$
the limit finite rotations $\Omega_k$ are mutually inverse. This property
gives a perfect discrete analog of limit permanent rotations in the
classical Suslov problem.

\medskip
\noindent{\it Proof of Theorem} \ref{discr_sus_dyn}. First, we
describe the discrete dynamics on the part ${\cal U}_0$ of the
momentum surface $\cal U$ bounded by the condition
$$
E=( J_{11}+J_{33}) M_1^2 +2 J_{12} M_1 M_2+ ( J_{22}+J_{33}) M_2^2 \le h^* .
$$
For this purpose introduce a new 
coordinate system
\begin{align*}
{\cal M}_1 &= ( J_{13}(J_{11}+J_{33})+J_{12}J_{23}) M_{1}+
(J_{23}(J_{22}+J_{33})+J_{12}J_{13}) M_{2}, \\
{\cal M}_2 &=J_{23} M_1 -J_{13}M_2 .
\end{align*}
In view of relations (\ref{full}) one has
\begin{align}
{\cal M}_1 &= ( J_{13}q_{1} + J_{23}q_{2}) \left[ Q +
((J_{11}+J_{33}) ( J_{22}+J_{33})- J_{12}^{2}) q_{0} \right ], \label{M_1}
\\
{\cal M}_2 &= Q q_0 -( J_{13}q_{1} + J_{23}q_{2})^2, \label{M_2} \\
Q &= (J_{12}J_{13}+J_{22}J_{23} +J_{23}J_{33}) q_{1}
-( J_{11}J_{13}+J_{12}J_{23}+J_{13}J_{33}) q_{2} . \nonumber
\end{align}
Using the properties $(J_{11}+J_{33}) ( J_{22}+J_{33})- J_{12}^{2}>0$,
$q_0\ge 0$, one can show that
in the domain ${\cal U}_0$ the expression in square brackets in (\ref{M_1})
is positive.
Hence, on the segment of the line ${\cal M}_1=0$ in ${\cal U}_0$ on has
$J_{13}q_{1} + J_{23}q_{2}=0$ and it consists of stationary points of the
map.
The points of ${\cal U}_0$ with positive (negative) ${\cal M}_1$ correspond
to the points
on the interior part of ${\cal S}\subset \mathbb{RP}^2$ with positive
(respectively negative)
values of $J_{13}q_{1} + J_{23}q_{2}$.
Next, in view of (\ref{diff_suslov}),
$$
{\cal M}_{2,k+1} = {\cal M}_{2,k} + ( J_{13}q_{1} + J_{23}q_{2})^2,
$$
which implies that the coordinate ${\cal M}_{2}$ always increases while
the point $M_k$ approaches the line ${\cal M}_1=0$ along the ellipse $E(M_1,
M_2)= \text{const} $.
Then, as follows from (\ref{M_2}), for $k\to -\infty$, one has ${\cal
M}_{2}<0$, $Q<0$
and for $k\to \infty$, ${\cal M}_{2}>0$, $Q>0$. As a consequence, the
equilibria positions
on ${\mathbb P}_u ={\mathbb P} \cap \Theta_{-}$ are unstable and those
on ${\mathbb P}_s ={\mathbb P} \cap \Theta_{+}$ are stable.

Further, due to (\ref{diff_suslov}),
${\cal M}_{1,k+1} - {\cal M}_{1,k} = -( J_{13}q_{1} + J_{23}q_{2}) Q$
and, therefore,
$$
{\cal M}_{1,k+1}
= ( J_{13}q_{1} +J_{23}q_{2}) ( (J_{11}+J_{33})(J_{22}+J_{33})- J_{12}^{2})
q_{0} .
$$

The latter and \eqref{M_1} implies that, unless $J_{13}q_{1}
+J_{23}q_{2}=0$, the coordinates
${\cal M}_{1,k}$ and ${\cal M}_{1,k+1}$ always have the same sign, i.e., the
sequence
$\{M_k\}$ lies entirely in one of the domains ${\cal U}_0\cap \{ {\cal
M}_1\lessgtr 0\}$.
Reformulating these properties for the interior part of the domain ${\cal
S}\subset \mathbb{RP}^2$,
we arrive at the statement of the theorem.

 
\section{Discrete Unbalanced Chaplygin Sleigh}
Now we pass to discretization of the EPS equations (\ref{EPS.se2.eqn}) 
on the coalgebra $se^*(2)$.

The two subsequent positions of the sleigh are given by the matrices
\[
X_{k}=\left( 
\begin{array}{ccc}
\cos \theta _{k} & -\sin \theta _{k} & x_{k} \\ 
\sin \theta _{k} & \cos \theta _{k} & y_{k} \\ 
0 & 0 & 1%
\end{array}%
\right) ,\qquad X_{k+1}=\left( 
\begin{array}{ccc}
\cos \theta _{k+1} & -\sin \theta _{k+1} & x_{k+1} \\ 
\sin \theta _{k+1} & \cos \theta _{k+1} & y_{k+1} \\ 
0 & 0 & 1%
\end{array}%
\right) 
\]
The {\it helical displacement} in the body frame is defined by 
$ \varOmega_k = X_{k}^{-1} X_{k+1} \in SE(2)$ and 
straightforward computation shows that
\begin{gather} 
\label{Om_k}
\varOmega_k = \begin{pmatrix}
\cos ( \Delta \theta_k) & - \sin ( \Delta \theta_k) &  
\cos \theta_k \,\Delta x _k + \sin \theta_k \,\Delta y _k \\
\sin ( \Delta \theta_k)  &\hphantom{-} \cos ( \Delta \theta_k) 
& - \sin \theta_k \,\Delta x _k + \cos \theta_k \,\Delta y _k \\
0 & \ \ 0 & \ \ \ 1
\end{pmatrix} , \\
\Delta \theta_k = \theta_{k+1}-\cos \theta_{k}, \quad 
\Delta x _k = x_{k+1}-x_{k}, \quad \Delta y_k = y_{k+1}-y_{k} .
\nonumber
\end{gather}

Following the expression (\ref{approx}), define the left-invariant discrete Lagrangian on
$SE(2)\times SE(2)$ by replacing the helical velocity $\xi$ in
(\ref{cont.sl.lagr}) with $X_{k}^{-1}(X_{k+1}-X_k)$.  Up to an additive constant, we get
\begin{gather} \label{sl}
L_d (X _{k+1}, X _k) = \frac12 \mbox{tr } 
\left(\varOmega_k {\mathbb J} \varOmega_k^T \right) -  \frac12 \mbox{tr } 
\left( {\mathbb J} \varOmega_k^T  + \varOmega_k {\mathbb I}\right), \\
{\mathbb J}= \begin{pmatrix}
J/2+ma^{2} & m ab & ma \\ 
m ab  & J2+mb^{2} & mb \\ 
ma & mb & m
\end{pmatrix} , \nonumber
\end{gather} 
where, as above, $a,b$ are coordinates of the mass center $C$ in the body frame and
$J$ is its scalar moment of inertia with respect to the origin $A$.
This yields the following scalar expression
\begin{align}
 L_d =&\frac{m}{2} \Delta y_k^{2}
+\frac{m}{2} \Delta x_k^{2}+\left( J+ma^{2}+mb^{2}\right) (1-\cos \Delta\theta_{k})\nonumber
 \\
&\quad +am[ ( \sin \theta_{k+1}-\sin \theta _{k} ) 
\Delta y_k +( \cos \theta _{k+1}-\cos \theta _{k}) \Delta x _k ] \nonumber \\
& \quad +bm[ ( \cos \theta _{k+1}-\cos \theta _{k}) \Delta y_k 
-( \sin \theta _{k+1}-\sin \theta _{k}) \Delta x _k ] . \label{discr.L2}
\end{align}
In the continuous limit, when 
\begin{gather} \label{cont.lim}
\Delta\theta _{k}= \varepsilon \, \omega+ O(\varepsilon^2), 
\quad \Delta x_{k} =\varepsilon \dot x + O(\varepsilon^2),  \quad
\Delta y_{k} =\varepsilon \dot y + O(\varepsilon^2),\qquad \varepsilon <<1 , \\
\cos \theta _{k+1}-\cos \theta _{k}=- \varepsilon\, \omega \sin \theta + O(\varepsilon^2), \quad
\sin \theta _{k+1}-\sin \theta _{k}= \varepsilon\, \omega \cos \theta + O(\varepsilon^2), 
\nonumber
\end{gather}
expression (\ref{discr.L2}) divided by $\varepsilon$ transforms to the continuous 
Lagrangian (\ref{Caplag}) plus higher order terms in $\varepsilon$. 
\medskip

According to definition (\ref{disc.mom.dfn}), the discrete momentum in the body \\
$P_k=( p_{\theta,k}, p_{1,k}, p_{2,k}) \in se^*(2)$, has the form
\begin{align*}
p_{\theta,k} & = - \pder {}{\varepsilon} L_d 
( {\theta_k}+ \varepsilon ,\theta_{k+1}, x_k, x_{k+1} ,y_k,y_{k+1} ) \bigg |_{\varepsilon=0 } \\
p_{1,k} &=  - \pder {}{\varepsilon}  L_d 
( {\theta_k} ,\theta_{k+1} , x_k+\varepsilon \cos\theta_{k}, x_{k+1} , 
y_k +\varepsilon \sin\theta_{k}, y_{k+1})\bigg |_{\varepsilon=0 }  , \\
p_{2,k} &= - \pder {}{\varepsilon}  L_d 
( {\theta_k}, \theta_{k+1}, x_k-\varepsilon \sin\theta_{k}, x_{k+1}, 
y_k+\varepsilon \cos\theta_{k},y_{k+1} )\bigg |_{\varepsilon=0} , 
\end{align*} 
that is, 
\begin{align}
p_{\theta,k } & = (J+m a^{2}+m b^{2})\sin (\Delta\theta_{k})+am V_{2,k} -bm V_{1,k}\, , 
\nonumber \\
p_{1,k} &= mV_{1,k} - am (1-\cos (\Delta\theta_{k}) )- bm \sin (\Delta\theta_{k})  ,
\label{disc.mom.sl} \\
p_{2,k} &= mV_{2,k} +am \sin (\Delta\theta_{k}) - bm  (1-\cos (\Delta \theta_{k}) ) \nonumber,
\end{align} 
where 
\begin{align}
V_{1,k} &= (\varOmega_k)_{13}\equiv \Delta x_k \cos\theta_{k}+\Delta y_k \sin \theta_{k} , \nonumber \\ 
V_{2,k} &= (\varOmega_k)_{23}\equiv - \Delta x_{k} \sin\theta_{k}+\Delta y_{k}\cos \theta_{k} 
\label{discr.vel}
\end{align} 
are ''discrete velocities'' in the body frame.
\medskip

Next, 
the coadjoint action on $se^*(2)$ can be written in form
\begin{equation} \label{discr.ad.sl} 
\mbox{Ad}_{\varOmega_k}^* P_{k} = 
\begin{pmatrix}
 &  p_{\omega,k} - p_{2,k} V_{1,k} + p_{1,k} V_{2,k}  \\
 & \cos (\Delta\theta_{k}) p_{1,k}  + \sin (\Delta\theta_{k})p_{2,k} \\
 & -\sin (\Delta\theta_{k}) p_{1,k}+ \cos (\Delta\theta_{k}) p_{2,k}  
\end{pmatrix} . 
\end{equation}

In the absence of constraints the dynamics of the 2-dimensional body can be 
represented by the discrete Euler--Poincar\'e equations
\begin{equation} \label {discrete.free.se2.eqn}
P_{k+1} = \mbox{Ad}_{\varOmega_k}^* P_{k},
\end{equation} 
which gives the momentum conservation law written in the body
frame. In particular, for $a=b=0$ (the mass center $C$ lies at the origin), 
the system (\ref{discrete.free.se2.eqn}),  (\ref{discr.ad.sl}) yields
\begin{align*}
\sin (\theta _{ k+1} - \theta _k) & = \sin (\theta _k - \theta _{k-1}),\\ 
\Delta x_{k+1} \cos\theta_{k+1}+\Delta y_{k+1} \sin \theta_{k+1}  
&= \Delta x_{k} \cos\theta_{k+1}+\Delta y_{k} \sin \theta_{k+1},\\
- \Delta x_{k+1} \sin\theta_{k+1}+\Delta y_{k+1}\cos \theta_{k+1}  
 &= -\Delta x_{k} \sin\theta_{k+1}+\Delta y_{k}\cos \theta_{k+1}  ,
\end{align*}
which implies that {\it for small $\theta$'s\/} the differences
$\theta _{k+1} - \theta_{k} $, $ x_{k+1} - x_{k}$, 
and $ y_{k+1} - y_{k}$ are the same for any integer $k$, the result one expects
from studying the continuous problem.

\paragraph{Discrete constraint on $SE(2)$.}
We now impose discrete left-invariant constraints on $SE(2)\times SE(2)$
in the form of restrictions on discrete helical velocities $\varOmega_k = X_{k+1} X_{k}^T$.
By analogy with continuous constraint defined by (\ref{se2.constr}), 
a naive choice of a discrete constraint is just to set
\begin{equation}
\label{false!}
(\varOmega_k)_{23}\equiv- \sin \theta _k \,\Delta x _k + \cos \theta _k
\,\Delta y _k = 0.
\end{equation}
This choice however is not the right one. 
Indeed, following our approach to discrete left-invariant constraints, 
admissible rotations and translations
must be exponents of the matrices of the form (\ref{d_for_se_n}). In
this case ${\mathfrak h}$ generates the subgroup $SE(n-1)$ and, according to 
Proposition \ref{S-Cartan}, 
$\exp{\mathfrak d}$ must be a covering of the homogeneous space
$SE(n)/SE(n-1)$. 

In the particular case $n=2$, when $\varOmega_k$ is given by (\ref{Om_k}), we have
\begin{proposition} \label{circ.motion}
The variety ${\mathcal S}=\exp{\mathfrak d}\subset SE(2)$ is diffeomorphic to the
the canonical line bundle $\pi\, : \, {\cal L}\to {\mathbb RP}^1=(z_1:z_2)$ (Moebius cylinder)
such that
$\pi^{-1} (z_1:z_2) =\{ \nu (z_1,z_2), \nu\in {\mathbb R}\}$ and it is defined by the condition
\begin{equation} \label{constr.omegas}
\frac{\varOmega_{23}}{\varOmega_{13}} =\frac{1-\varOmega_{11}}{\varOmega_{21}} \, .
\end{equation}
The latter yields the following constraint
\begin{align}
- (\Delta x _k \cos \theta _k + \Delta y _k & \sin \theta _k)
\sin (\Delta \theta_k/2) \nonumber \\
 & \quad + ( - \Delta x _k \sin \theta _k + \Delta y _k \cos \theta _k)
\cos(\Delta \theta_k/2) = 0 \label{se2.const}
\end{align}
or, equivalently, 
\begin{equation} \label{se.const}
V_{1,k}[1-\cos(\Delta\theta_k) ] - V_{2,k} \sin(\Delta\theta_k)= 0 .
\end{equation}
\end{proposition}

The corresponding left-invariant constraint on $SE(2)\times SE(2)$ has the form
\begin{equation} 
\label{discrete.se2.constraint.eqn}
- \sin \left(\frac{\theta _{ k+1} + \theta_k}{2}\right)
(x_{k+1} - x _k ) + \cos \left(\frac{\theta_{k+1}
+ \theta_k}{2}\right) (y_{k+1} - y _k )=0 .
\end{equation}
Observe that in the continuous limit (\ref{cont.lim}) this yields the
constraint $-\dot x\sin\theta+\dot y\cos\theta=0$. 
\medskip

\noindent{\it Proof of Proposition} \ref{circ.motion}
For an element $S\in {\mathfrak d}$ we have
$$
S=\left(
\begin{array}{ccc}
0 & -\omega & v \\
\omega & 0 & 0 \\
0 & 0 & 0%
\end{array}\right), \quad
\exp (S t)=\left(
\begin{array}{ccc}
\cos \omega t & -\sin \omega t & \frac{v}{\omega } \sin (\omega t) \\
\sin \omega t & \cos \omega t & \frac{v}{\omega }( 1-\cos (\omega t) ) \\
0 & 0 & 1
\end{array} \right) ,
$$
where $\omega, v$ are arbitrary.
As a result, for the points of the admissible shift subvariety, relation (\ref{constr.omegas})
holds. Next, since
\begin{equation} 
\label{quotient}
\omega t= \Delta\theta_k, \quad \mbox{and} \quad
\frac{\varOmega_{23}}{\varOmega_{13}}\equiv \frac{1-\cos\omega t} {\sin\omega t}
= \tan\frac{\Delta\theta_k} 2 ,
\end{equation} 
in view of (\ref{discr.vel}), we have (\ref{se2.const}) and (\ref{se.const}). 

Finally, as seen from the last relation, the angle $\Delta\theta_k$ determines
the quotient $\varOmega_{23}/\varOmega_{13}$, i.e., a line  in
${\mathbb R}^2=(\varOmega_{23}, \varOmega_{13})$. As $\Delta\theta_k$ changes by $2\pi$,
the line rotates by $\pi$, hence $\cal S$ is diffeomorphic to the Moebius cylinder.
\medskip

\begin{figure}[h,t]
\begin{center}
\includegraphics[width=0.4\textwidth]{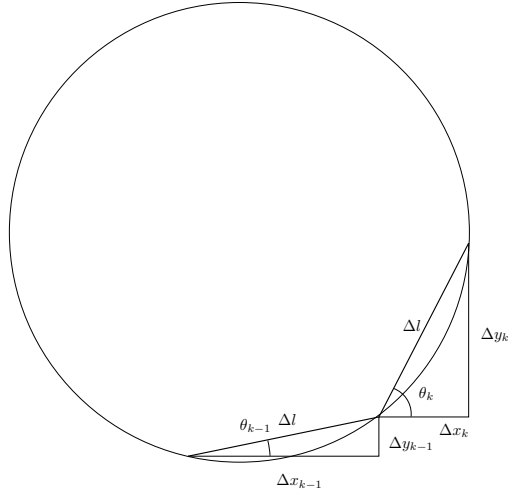}
\caption{The geometry of the incremental displacements for the
  Chaplygin sleigh.}
\label{circ.fig}
\end{center}
\end{figure}

\paragraph{Remark.}
As seen from relation (\ref{se2.const}), the matrices from ${\cal
  S}\subset SE(2)$ describe ``circular translations'' of the sleigh along
the axis $X_k$ of the blade: the points $(x_k,y_k)$ and
$(x_{k+1},y_{k+1})$ in ${\mathbb R}^2$ must lie on a circle
such that the lines $X_k$ and $X_{k+1}$ are tangent to this
circle. This property also implies that
\begin{equation} 
\begin{aligned}
\Delta x_k \cos\theta_{k}+\Delta y_k \sin \theta_{k}
 &= \Delta x_k \cos\theta_{k+1}+\Delta y_k \sin \theta_{k+1},  \\
 -\Delta x_{k} \sin\theta_{k} +\Delta y_{k}\cos \theta_{k}
 &=  \Delta x_{k} \sin\theta_{k+1}- \Delta y_{k}\cos \theta_{k+1}
\end{aligned}
\label{discr.vel1}
\end{equation} 
(see Figure \ref{circ.fig}).  

The above constraint has also the following interpretation: in order
to transfer the sleigh from $(\theta_k, x_k, y_k)\in SE(2)$ to
$(\theta_{k+1}, x_{k+1}, y_{k+1}) \in SE(2) $ (assuming that this
transition is possible), one needs first to perform the rotation over
$\Delta\theta_k/2$ at $(x_k,y_k)$, which aims the sleigh towards
$(x_{k+1}, y_{k+1})$, then slide the sleigh from $(x_k,y_k)$ to
$(x_{k+1}, y_{k+1})$, and then perform another rotation over
$\Delta\theta_k/2$ (now at $(x_{k+1},y_{k+1})$).
\medskip

Note that under the constraint (\ref{se.const}) the image of the
discrete Legendre transformation (\ref{disc.mom.sl}) is an algebraic
quartic subvariety $\cal U$ in $se^*(2)=( p_{\theta}, p_{1}, p_{2})$
and to a generic pair $(p_{\theta}, p_{1})$ there correspond four
distinct points on $\cal U$ and four inverse images on ${\cal
  S}\subset SE(2)$.

\paragraph{Discrete momentum locus ${\cal U}\subset se^*(2)$.}
Below we concentrate on the important case $b=0$, when the structure of the
real surface ${\cal U}\subset se^*(2)$ becomes simpler. It is more
convenient to describe the image $\tilde {\cal U}$ of ${\cal U}$
in ${\mathbb R}^3= (p_\theta, \hat p_1, z)$, where 
$\hat p_1=a p_1 +2ma^2$, $z= \sin (\Delta\theta)$.

\begin{lemma} \label{surface}
\begin{description}
\item{1).} For $b=0$ the surface $\tilde{\cal U}$ is given by cubic equation
\begin{equation} \label{det.delta.phi}
{\cal H}(p_\theta,\hat p_1,z)= J^2 z^3 -2 J p_{\theta,k} z^2
 +( \hat p_1^2 + 2J \hat p_1  + p_{\theta}^2)z - 2p_{\theta}\hat p_1=0 .
\end{equation}
$\tilde{\cal U}$  lies entirely between the planes $z=\pm 1$ and is tangent to them along the
lines $\ell_\pm=\{\pm p_\theta-\hat p_1= J \}$ respectively. The $p_\theta$- and 
$\hat p_1$-axis belong entirely to $\tilde{\cal U}$.

\item{2).} For the parts of $\tilde{\cal U}$ over the quadrants
\begin{align*}
L_{++} & = \{- \hat p_1  + p_\theta > J\} \cap \{ -\hat p_1  - p_\theta > J\} \quad
\textup{and} \\
L_{--} &= \{-\hat p_1 + p_\theta < J\} \cap \{-\hat p_1- p_\theta < J\}
\end{align*}
one has $\cos(\Delta\theta)>0$, i.e.,
$-\pi/2 <\Delta\theta<\pi/2$ and in the rest of quadrants one has $\cos(\Delta\theta)<0$
$(\pi/2<\Delta\theta<3\pi /2)$.


\item{3).}
The projection $\Pi \, : \, {\tilde U}\to {\mathbb R}^2=(p_\theta,
\hat p_1)$ is one-to-one except the above segments and the interior of
triangular domain bounded by the discriminant curve
$$
\hat p_1^4 + 6J \hat p_1^3 +\hat p_1^2(12 J^2+2p_\theta^2) -\hat p_1 (10 J p_\theta^2 -8J^3)
+p_\theta^4-J^2 p_\theta^2 =0.
$$
The curve is symmetric with respect to $\hat p_1$-axis, it is tangent to
$p_\theta$-axis at the origin $(0,0)$ and has 3 cusp points with
coordinates $(0,-2J)$, $(c_1, c_2), (-c_1, c_2)$, with some
positive constants $c_1, c_2$. In this domain the projection $\Pi$ is 3 to 1.

\item{4).} The curve $\{V_1=0\}\subset \tilde{\cal U}$ is projected onto the ellipse 
\begin{equation} 
\label{ell}
{\cal E}= \{ p_\theta=(J+ma^2) \sin (\Delta\theta), \quad 
\hat p_1 = ma^2(1+\cos (\Delta\theta)) \mid \Delta\theta\in (0;2\pi)\}. 
\end{equation}
Inside the ellipse the values of $V_1$ are negative and outside are positive.
\end{description}
\end{lemma}

\begin{figure}[h,t]
\begin{center}
\includegraphics[width=0.45\textwidth]{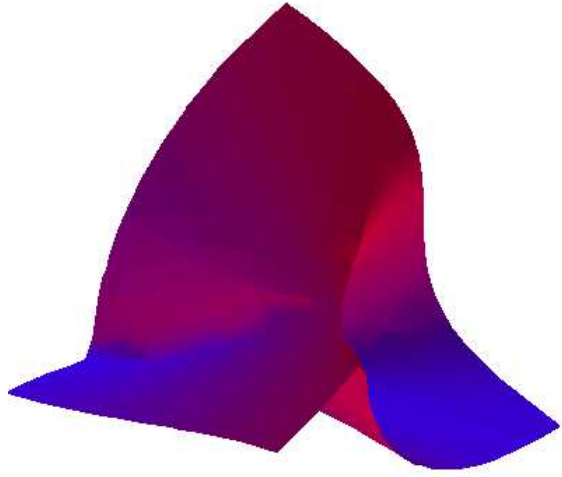}
\caption{The surface $\tilde{ \cal U}$.}
\label{cubic.fig}
\end{center}
\end{figure}

Note that the point $O$ with coordinates $p_\theta=0$, $\hat
p_1=2ma^2$ corresponds to the origin in the $(p_\theta, p_1)$ phase
plane and in a neighborhood of this point the projection $\Pi$ is
one-to-one.  An example of the surface $\tilde{\cal U}$
for 
$J=1.5$ is presented in Figure
\ref{cubic.fig}.

%

\noindent{\it Proof of Lemma} \ref{surface}. 1). Using the condition
(\ref{se.const}), we exclude $V_1, V_2$ from the first two equations
of (\ref{disc.mom.sl}) to obtain the following condition on
$\Delta\theta_{k}$,
\begin{equation} \label{quartic}
J \sin^2 (\Delta\theta_{k}) - p_{\theta,k} \sin(\Delta\theta_{k})
+ (2ma^2+a p_{1,k})[1- \cos (\Delta\theta_{k}) ] =0.
\end{equation}
This equation always has trivial solution $\Delta\theta_{k}=2n \pi$,
$n\in {\mathbb Z}$.  Setting $z= \sin(\Delta\theta_{k})$, $\cos
(\Delta\theta_{k})=\sqrt{1-z^2}$, $\hat p_1=a p_x + 2ma^2$, we arrive
at a quartic polynomial equation with respect to $z$, which has the
root $z=0$. Factoring it out and omitting the index $k$, one gets the
cubic equation (\ref{det.delta.phi}).

Now setting in (\ref{det.delta.phi}) $z=\pm 1$, we get $(J \mp p_\theta +\hat p_1)^2=0$, which
implies that
$\tilde{\cal U}$ is indeed tangent to the planes $z=\pm 1$  along the lines $\ell_\pm$.
Finally, setting $z=p_\theta=0$ or  $z=\hat p_1=0$, one sees that equation
(\ref{det.delta.phi}) is satisfied for any $\hat p_1$ and $p_\theta$ respectively.

2). For fixed $p_\theta, \hat p_1$, each root of (\ref{det.delta.phi})
gives a solution of (\ref{quartic}) with a sign of $\cos(\Delta\theta_{k})$ appropriately chosen.
As seen from (\ref{quartic}), for large $|p_x|$ and small $z=\sin(\Delta\theta)$, the
value of $\cos (\Delta\theta)$ must be close to 1,
whereas for large $p_\theta>0$ and small $|p_x|$,
$\cos (\Delta\theta)$ must be negative. Since the sign of $\cos (\Delta\theta)$ can
change only under passage from one quadrant on the plane $(p_\theta, Y)$ to another one,
this proves item 2).

Items 3), 4) are verified by direct calculations.
\medskip

\paragraph{Discrete dynamics on $se^*(2)$ with the constraint.}
According to (\ref{disc.momentum.eqn}), 
the discrete Euler--Poincar\'e--Suslov equations associated with the constraint
(\ref{se2.constr}) have the form
\begin{equation} \label{EP_sleigh}
P_{k+1} = \mbox{Ad}_{\varOmega_k}^* P_{k}+ \lambda_k (0,0,1)^T .
\end{equation}

Substituting here expressions (\ref{discr.ad.sl}), we find that 
under the constraint (\ref{se.const}) the first two components of $P_{k+1}$ have the form
\begin{align*}
p_{\theta,k+1} & = (J+m a^{2}+m b^{2})\sin (\Delta\theta_{k})-bm V_{1,k} \\
& \quad + am[- \Delta x_{k} \sin\theta_{k+1}+ \Delta y_{k}\cos \theta_{k+1}] , \\
p_{1,k+1} & =mV_{x,k} - bm \sin (\Delta\theta_{k}) \\
 & \quad + am[\Delta x_{k} \cos\theta_{k+1}+ \Delta y_{k}\sin \theta_{k+1}],  
\end{align*}
which, in view of (\ref{discr.vel1}), (\ref{disc.mom.sl}), yields 
\begin{equation}
\label{discrete.constr.se2.eqn} 
\begin{aligned}
p_{\theta,k+1} & = p_{\theta,k}-2 am V_{2,k}\, ,  \\
p_{1,k+1} &= p_{1,k} + 2am(1-\cos (\Delta\theta_{k}))\, . 
\end{aligned}
\end{equation}
Expressions (\ref{discrete.constr.se2.eqn}), (\ref{se.const}) define a multi-valued map
$\cal U \to \cal U$ or $\cal S \to \cal S$ which consists of 3 steps:
\begin{description}
\item{1).} Given $\Delta\theta_{k}, V_{1,k}$, one finds $ V_{2,k}$ form the constraint 
(\ref{se.const})
and then $p_{\theta,k}, p_{1,k}, p_{2,k}$ from the Legendre transformation (\ref{disc.mom.sl}).

\item {2).} One finds $p_{\theta,k+1}, p_{2,k+1}$ from (\ref{discrete.constr.se2.eqn}).

\item{3).} One finds $\Delta\theta_{k+1}, V_{1,k+1}$ 
by choosing a solution of the system of equations
\begin{align*}
p_{\theta,k+1} & = (J+m a^{2}+m b^{2}) \sin (\Delta\theta_{k+1})
+\left( am \frac{1-\cos(\Delta\theta_{k+1})}{\sin(\Delta\theta_{k+1})} -bm\right) V_{x,k+1}\, , \\
p_{1,k+1} &= mV_{x,k+1} - am (1-\cos (\Delta\theta_{k+1}) )- bm \sin (\Delta\theta_{k+1}) ,
\end{align*}
which are obtained from (\ref{disc.mom.sl}), (\ref{se.const}) by replacing $k\to k+1$. 
\end{description}
\medskip

\begin{theorem} \label{int_chapl}
Equations (\ref{discrete.constr.se2.eqn}) preserve the quantity 
\begin{equation} \label{dis.chapl.int}
E= m p_\theta^2 +2 bm p_\theta p_1 + (J + m(a^2 + b^2)) p_1^2 ,
\end{equation}
which coincides with the truncated energy integral (\ref{cont_sl_int}) 
of the continuous Chaplygin sleigh.
\end{theorem}

\noindent{\it Proof.} Substituting expressions (\ref{disc.mom.sl}) and
(\ref{discrete.constr.se2.eqn}) into (\ref{dis.chapl.int}) and taking
into account the constraint (\ref{se.const}), one obtains the same
expression in terms of $\Delta\theta_{k}, V_{1,k}$ and $V_{2,k}$.
\medskip

Since the quadratic form (\ref{dis.chapl.int}) is positive definite,
we conclude that the invariant manifolds of the map
\eqref{discrete.constr.se2.eqn} are the ellipses in the $ p_\theta p_1
$-plane.

\paragraph{Stationary solutions of the discrete Chaplygin sleigh.} 
As follows from (\ref{discrete.constr.se2.eqn}), 
for the initial conditions $\{\Delta\theta_k=0, V_{2,k}=0\}$ one has
$$
p_{\theta,k+1}=p_{\theta,k},\quad  p_{1,k+1}=p_{1,k}.
$$
Hence, it is natural to choose  such $\lambda_k$ in (\ref{EP_sleigh}) that 
$p_{2,k+1}= p_{2,k}$ as well.
Thus, like the continuous system (\ref{ch_sl_1.eqn}), for $a\ne 0$
the map (\ref{EP_sleigh}) has a family of 
stationary solutions which, on the momenta plane 
$(p_\theta, p_1)$, is represented by the line $\{p_\theta+bp_1=0 \}$. 
Such solutions correspond to shifts in the $(x,y)$-plane along the axis of
the blade by constant distances. 

On the other hand, for $a=0$ all the solutions are stationary.
That is, in contrast to the case of absence of constraints, 
when the discrete momentum in space is preserved,
now it is the momentum in the body $P$, which is preserved. 
In view of (\ref{disc.mom.sl}), this implies
$$
\Delta\theta_{k+1}= \Delta\theta_{k}, \quad V_{1,k+1}=V_{1,k} .
$$
As a result, the discrete trajectory on the plane $(x,y)$ consists of
displacements along a circle
with radius $\rho=V_{1,k}/\sin(\Delta\theta_{k})$
\footnote{ As numerical simulations show, if one chooses the naive constraint 
(\ref{false!}) instead of (\ref{constr.omegas}), 
then for $a=0$ the discrete trajectory on the plane $(x,y)$ lies on a spiral.}
. The same behavior occurs to the continuous sleigh for $a=0$.

\paragraph{The case $b=0$, $a\ne 0$.}
In this case the map $(p_{\theta,k}, p_{x,k}) \to (p_{\theta,k+1}, p_{1,k+1})$
has a line of stationary points $p_\theta=0$, and,
according to Theorem \ref{int_chapl}, the discrete trajectories lie on symmetric
invariant ellipses $ m p_\theta^2 + (J + ma^2) p_1^2 =E$. 
Without loss of generality, we assume $a>0$. Then the following property holds.

\begin{theorem} In the neighborhood of the origin $O$ bounded by the condition 
$E < m ^2 a ^2 (J + ma^2) $ the map is single-valued and has a
  bi-asymptotic behavior similar to that of the continuous Chaplygin
  sleigh system. Namely, for $k\to -\infty$, the point $(p_{\theta,k},
  p_{x,k})$ approaches, along the corresponding invariant ellipse, a
  point of the 
  segment $\{p_\theta=0,\; -ma <p_1 < 0\}$ of unstable stationary
  points, and for $k\to +\infty$, the point $(p_{\theta,k},
  p_{x,k})$ approaches one of the points of the segment
  $\{p_\theta=0,\; 0<p_1 <ma\}$ of stable stationary points. In both
  cases the sequence $\{ (p_{\theta,k}, p_{x,k})\}$ remains in one of
  the half-planes $p_\theta<0$ or $p_\theta>0$.
\end{theorem}

\noindent{\it Proof.} Part (3) of Lemma \ref{surface} implies that the
map is single-valued in the region $ E < m ^2  a ^2 (J + ma^2)  $.
 
Next, as follows from the first relation in
(\ref{discrete.constr.se2.eqn}) for $a>0$, the increment
$p_{1,k+1}-p_{1,k}$ is always greater than or equal to zero.  Then, to
prove the bi-asymptotic behavior, it remains only to show that the
sequence $\{(p_{\theta,k},p_{1,k})\}$ lies entirely in one of the
half-planes $p_\theta \lessgtr 0$.

Indeed, let the point $(p_{\theta,k}, p_{x,k})$ be inside the ellipse
$\cal E$ given by (\ref{ell}). First, assume that $p_{\theta,k}>0$.
Then, in view of items (2), (4) of Lemma \ref{surface}, and the
constraint (\ref{se2.const}), $V_{1,k}$ and $V_{2,k}$ are
negative. According to (\ref{discrete.constr.se2.eqn}), the increment
$p_{\theta,k+1}- p_{\theta,k}$ is then positive.  Similarly, for
$p_{\theta,k}<0$ one has $p_{\theta,k+1}- p_{\theta,k}<0$.

Next, if $0 < p_{1,k}<ma$ and $(p_{\theta,k}, p_{x,k})$ lies in the domain 
$ E<(J + ma^2) m^2 a^2$, then, using (\ref{disc.mom.sl}), (\ref{se2.const}),
one shows that $ 2 am V_{2,k} > p_{\theta,k} $ for $ p_{\theta,k}>0 $
and $ 2 am V_{2,k} < p_{\theta,k} $ for $ p_{\theta,k}<0 $.
Therefore, in view of (\ref{EPS.se2.eqn}), 
$p_{\theta,k+1}>0$, respectively, $p_{\theta,k+1}<0$. 

As a result, in any case, $p_{\theta,k}$ and $ p_{\theta,k+1}$ cannot
have different signs, which proves the theorem.
\medskip

One concludes that in the neighborhood of the origin $O$ the
discrete-time dynamics is similar to that of the Suslov problem
illustrated in Figure \ref{quat_dots}.

\begin{figure}[h,t]
\begin{center}
\includegraphics[width=0.6\textwidth]{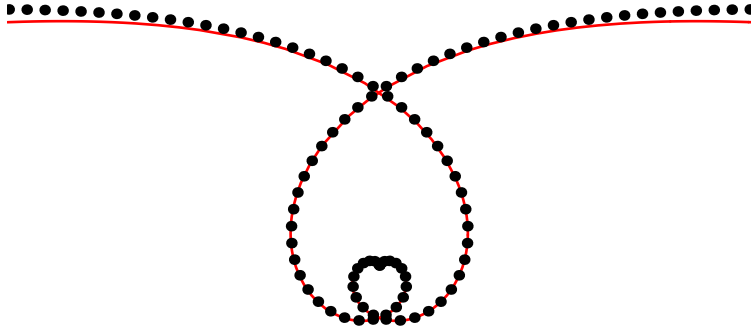}
\caption{A typical discrete sleigh trajectory.}
\label{cusps.fig}
\end{center}
\end{figure}

We conclude this section with an example of the discrete sleigh trajectory
on the plane $(x,y)$ compared to a continuous trajectory for $b=0$
with a cusp, as presented in Figure \ref{cusps.fig}.

\section{Conclusions} 
The discrete nonholonomic Suslov problem and the Chaplygin sleigh that
we introduced in this paper 
properties of their corresponding continuous-time dynamical systems;
in particular, they preserve the reduced constrained energy and, in
the balanced case, the momentum.  This enables one to obtain explicit
solutions for the momentum dynamics of both discrete systems in terms
of theta-functions and exponents.
It is not currently clear if the complete solvability is due to the
(low) dimension of the systems and if it is possible to construct
completely solvable discretizations of the multidimensional Suslov and
Chaplygin problems.  These issues will be addressed in a future
publication.

On the other hand, by modifying our approach, one can also consider
discretizations of nonholonomic LR systems on Lie groups. For such
systems, the Lagrangian is left-invariant while the constraint
distribution is {\em right}-invariant. The discrete dynamics of such
systems, as well as the existence of their invariant measure, is
currently being developed and will be exposed in a future publication.
 
\subsection*{Acknowledgments}
YNF's research was partially
supported by Spanish Ministry of Science and Technology grant BFM
2003-09504-C02-02; DVZ's research was partially supported by NSF grant
DMS-0306017.


\begin{thebibliography}{99}

\bibitem{AKN}  Arnold, V.\,I., V.\,V. Kozlov, and A.\,I. Neishtadt [1989],
Mathematical Aspects of Classical and Celestial Mechanics. {\em
  Dynamical System III}, Springer-Verlag, New York.

\bibitem{Bloch2003}
Bloch, A.\,M. [2003], {\em Nonholonomic Mechanics and Control},
Interdisciplinary Applied Mathematics {\bf 24},
Springer-Verlag. 

\bibitem{BKMM} {Bloch A.\,M., P.\,S. Krishnaprasad, J.\,E. Marsden, and
  R. Murray} [1996], Nonholonomic Mechanical Systems with Symmetry.
  {\em Arch.\ Rational Mech.\ Anal.} \textbf{136}, 21--99.

\bibitem{BS1} Bobenko A.\,I. and Y.\,B. Suris [1999], Discrete Lagrangian
  Reduction, Discrete Euler--Poincar Equations, and Semidirect
  Products. {\em Lett.\ Math.\ Phys.} {\bf 49}, 79--93.


\bibitem{Ch1911} Chaplygin, S.\,A. [1911],
On the Theory of Motion of Nonholonomic Systems. The Theorem on the
Reducing Multiplier, {\em Math.\ Sbornik}
{\bf XXVIII}, 303--314, (in Russian).

\bibitem{CM1} Cort\'es J. and Mart\'inez S. [2001], Nonholonomic Integrators.
{\it Nonlinearity} {\bf 14}, 1365--1392.

\bibitem{FeKo2} {Fedorov Yu.\,N. and V.\,V. Kozlov} [1995],
Various Aspects of $n$-Dimensional Rigid Body Dynamics.
{\em Amer. Math. Soc. Transl} \textbf{168}, 141--171.

\bibitem{Gray} Gray, A. [1997], {\em Modern Differential Geometry of Curves and
Surfaces with Mathematica,} 2nd ed., Boca Raton, FL: CRC Press.

\bibitem{Jo2} Jovanovi\' c, B. [2001], Geometry and Integrability of
  Euler--Poincar\'e--Suslov Equations.  {\em Nonlinearity} \textbf{14}, 1555--1657.

\bibitem{Helgason} Helgason, S. [1962] Differential Geometry and Symmetric Spaces. 
Academic Press, New York.

\bibitem{Koz1} Kozlov, V.\,V. [1985],
On the Integration Theory of Equations of Nonholonomic
Mechanics. {\em Advances in Mechanics} {\bf 8}, 85--107 (in Russian).

\bibitem{Koz2} Kozlov, V.\,V. [1988] Invariant Measures of the
  Euler--Poincar\' e Equations on Lie algebras. {\em
    Funct.\ Anal.\ Appl.} {\bf 22}, 58--59.

\bibitem{deLeon1} de Le\'on, M., D. Mart\'in de Diego, and
  A. Santamar\'ia Merino. [2002]  Geometric Integrators and Nonholonomic
  Mechanics. arXiv:math-ph/0211028.

\bibitem{MCL1993}
McLachlan, R. [1993], Explicit Lie--Poisson Integration and the Euler
Equations, {\em Phys. Rev. Lett.} {\bf 71}, 3043--3046.

\bibitem{MarPekShk} Marsden, J.\,E., S. Pekarsky, and S. Shkoller [1999],
Discrete Euler--Poincar\'e and Lie--Poisson Equations.
{\em Nonlinearity} {\bf 12}, 1647--1662.

\bibitem{MaRa1999}
Marsden, J\,E. and T.\,S. Ratiu [1999], {\em Introduction to Mechanics
  and Symmetry,} Texts in Applied Mathematics {\bf 17}, Springer-Verlag.

\bibitem{MaWe2001}
Marsden, J.\,E. and M. West [2001], Discrete mechanics and variational
integrators, {\em Acta Numerica}, 357--514.

\bibitem{MosVes} Moser, J. and  A. Veselov [1991], Discrete Versions of Some
  Classical Integrable Systems and Factorization of Matrix
  Polynomials. {\em Comm.\ Math.\ Phys.} {\bf 139}, 217--243.

\bibitem{NeFu} {Neimark, Ju.\,I. and N.\,A. Fufaev} [1972] {\em Dynamics
  of Nonholonomic Systems}. Translations of Mathematical Monographs
  \textbf{33}, AMS, Providence.

\bibitem{No} {\it Steiner's Roman Surface.}
http://mathworld.wolfram.com/RomanSurface.html 

\bibitem{Su} {Suslov, G.} {\em Theoretical Mechanic}, Vol. 2, Kiev (in
  Russian).

\bibitem{Ve1988}
Veselov, A.\,P.,  [1988], Integrable Discrete-Time Systems and
Difference Operators, {\em Funk. Anal. Appl.} {\bf 22}, 1--13.

\bibitem{Ve1991}
Veselov, A.\,P.,  [1991], Integrable Lagrangian Correspondences and the
Factorization of Matrix Polynomials, {\em Funk. Anal. Appl.} {\bf 25},
38--49.

\bibitem{WeMa1997}
Wendland, J.\,M, and J.\,E. Marsden [1997], Mechanical Integrators
Derived from a Discrete Variational Principle, {\em Physica D} {\bf
  106}, 223--246.

\bibitem{Whitt} Whittaker, E.\,T. [1960] {\em A Treatise on Analytical
  Dynamics,} 4th ed., Cambridge Univ.\ Press, Cambridge.

\bibitem{Bl_Z}  Zenkov, D.\,V. and A.\,M. Bloch [2000],
Dynamics of the $n$-Dimensional Suslov problem. 
{\em J.\ Geom.\ Phys.} {\bf 34}, 121--136.

\bibitem{ZBl2002} Zenkov, D.\,V. and A.\,M. Bloch [2003], Invariant
  Measures of Nonholonomic Flows with Internal Degrees of
  Freedom. {\em Nonlinearity}\ {\bf 16}, 1793--1807.

\bibitem{Z2003}
Zenkov, D.\,V. [2003], Linear Conservation Laws of Nonholonomic
  Systems with Symmetry. {\em Discrete and Continuous Dynamical Systems
  (extended volume)}, 963--972.
\end{thebibliography}
\end{document}